\documentclass[10pt]{article}
\usepackage{amsfonts,color}
\usepackage{amsmath}
 \usepackage{epsfig}
\usepackage{lscape}
\usepackage{amssymb}
\usepackage{graphicx}
\usepackage{makeidx}
\usepackage{amssymb}
 \usepackage{footnote}

\DeclareMathOperator{\curl}{curl}

\DeclareMathOperator{\Curl}{Curl}

\newcommand{\Lp}[1]{L^{#1}(\Omega)}

\newcommand{\Hmp}[2]{H^{#1,#2}(\Omega)}
\newcommand{\Co}{C_0^{\infty}(\Omega)}

\newcommand{\norm}[1]{\|#1\|}
\newcommand{\R}{\mathbb{R}}

\newcommand{\C}{\mathbb{C}}

\renewcommand{\skew}{\mathop{\rm skew}}

\DeclareMathOperator{\sym}{sym}

\DeclareMathOperator{\dev}{dev}
\DeclareMathOperator{\sL}{\mathfrak{sl}}
\DeclareMathOperator{\so}{\mathfrak{so}}
\DeclareMathOperator{\gl}{\mathfrak{gl}}

\newcommand{\Sym}{ {\rm{Sym}} }

\newcommand{\Mprod}[2]{ {\langle #1 ,#2\rangle} }
\newcommand{\id}{ {1\!\!\!\:1 } }

\newcommand{\tr}[1]{ {\Tr \left[{#1}\right]} }

\allowdisplaybreaks[1]

\makeindex

\newtheorem{theorem}{Theorem}[section]

\newtheorem{corollary}[theorem]{Corollary}

\newtheorem{lemma}[theorem]{Lemma}

\def\barr{\begin{array}}
\def\id{1\!\!1}
\def\tr{\textrm{tr}}

\def\dvg{\textrm{Div}}
\def\crl{\textrm{Curl}}

\def\dd{\displaystyle}

\def\barr{\begin{array}}
\def\earr{\end{array}}
\def\bec#1{\begin{equation}\label{#1}}
\def\becn{\begin{equation*}}
\def\endec{\end{equation}}
\def\endecn{\end{equation*}}
 \def\C{\mathbb{C}}
 \def\H{\mathbb{H}}
  \def\L{\mathbb{L}_c}
  \def\cdp{\,\!\cdot\!\,}

\newtheorem{definition}{Definition}[section]
\newtheorem{remark}{Remark}[section]
\begin{document}

\title{The relaxed linear  micromorphic continuum:
existence, uniqueness and continuous dependence in dynamics}
\author{
Ionel-Dumitrel Ghiba\footnote{Ionel-Dumitrel Ghiba, \ \ \ \ Lehrstuhl f\"{u}r Nichtlineare Analysis und Modellierung, Fakult\"{a}t f\"{u}r Mathematik, Universit\"{a}t Duisburg-Essen, Campus Essen, Thea-Leymann Str. 9, 45127 Essen, Germany;  Alexandru Ioan Cuza University of Ia\c si, Department of Mathematics,  Blvd. Carol I, no. 11, 700506 Ia\c si,
Romania;  Octav Mayer Institute of Mathematics of the
Romanian Academy, Ia\c si Branch,  700505 Ia\c si; and Institute of Solid Mechanics, Romanian Academy, 010141 Bucharest, Romania, email: dumitrel.ghiba@uaic.ro} \quad
and \quad
Patrizio Neff\thanks{Patrizio Neff, \ \  Head of Lehrstuhl f\"{u}r Nichtlineare Analysis und Modellierung, Fakult\"{a}t f\"{u}r Mathematik, Universit\"{a}t Duisburg-Essen, Campus Essen, Thea-Leymann Str. 9, 45127 Essen, Germany, email: patrizio.neff@uni-due.de}
  \quad
and \quad\\
Angela Madeo\footnote{Angela Madeo, \ \  Laboratoire de G\'{e}nie Civil et Ing\'{e}nierie Environnementale,
Universit\'{e} de Lyon-INSA, B\^{a}timent Coulomb, 69621 Villeurbanne
Cedex, France; and International Center M\&MOCS ``Mathematics and Mechanics
of Complex Systems", Palazzo Caetani,
Cisterna di Latina, Italy,
 email: angela.madeo@insa-lyon.fr}
\quad
and \quad
Luca Placidi\footnote{Luca Placidi,\ \  Universit\`{a} Telematica Internazionale Uninettuno, Corso V. Emanuele II
39, 00186 Roma, Italy; and
International Center M\&MOCS ``Mathematics and Mechanics of Complex
Systems", Palazzo Caetani, Cisterna di Latina, Italy,
 email: luca.placidi@uninettunouniversity.net}
\quad
and \quad
Giuseppe Rosi\footnote{Giussepe Rosi,\ \  Laboratoire Mod\'{e}lisation Multi-Echelle, MSME UMR 8208 CNRS, Universit\'{e}
Paris-Est, 61 Avenue du General de Gaulle,
Creteil Cedex 94010, France; and  International Center M\&MOCS, University of
L'Aquila, Palazzo Caetani, Cisterna di Latina, Italy, email:
giuseppe.rosi@uniroma1.it}  }
\date{\it Dedicated to Prof. Antonio Di Carlo in recognition of his academic
activity}
\maketitle

\begin{abstract}
We study well-posedness  for the  relaxed linear elastic micromorphic continuum model with symmetric Cauchy force-stresses and curvature contribution depending only on the micro-dislocation tensor. In contrast to classical micromorphic models our
free energy  is not uniformly pointwise positive definite in the
control of the independent constitutive variables.  Another interesting feature concerns the prescription  of boundary values for the micro-distortion field: only tangential traces may be determined which are weaker than the usual strong anchoring boundary condition. There, decisive use is made of  new coercive
inequalities recently proved by Neff, Pauly and  Witsch and by Bauer, Neff, Pauly and Starke.   The new relaxed micromorphic formulation can be related to dislocation dynamics, gradient plasticity and seismic processes of earthquakes.
\\
\vspace*{0.25cm}
\end{abstract}

\newpage

\tableofcontents

\section{Introduction}

In this paper we show the well-posedness of a recently introduced new variant of the micromorphic model \cite{NeffGhibaMicroModel}. Micromorphic elasticity \cite{Eringen99,Neff_STAMM04,Mariano05,Mariano08a,Mariano08b,SteigmannIJNM12}
 is a generalized continuum formulation which tries to incorporate {\it microstructure} into the formulation of elasticity problems. This is necessary if one wants to describe size-effects (smaller is relatively stiffer), dispersion of waves phenomena etc. One of the best known such extension is the Cosserat model \cite{Cosserat09,Neff_ZAMM05,Neff_JeongMMS08,Neff_Jeong_Conformal_ZAMM08, Neff_Paris_Maugin09}. Our micromorphic  model equations are linear and the question is permitted as to what new kind of model there can be after the general framework has been introduced by Mindlin and Eringen \cite{Mindlin64,Eringen64,Eringen99}. Indeed, our new relaxed micromorphic model is a subclass of the classical model which, however, violates pointwise uniform definiteness of the energy: the new energy is positive semi-definite only.

The relaxed micromorphic model \cite{NeffGhibaMicroModel}  preserves full kinematical freedom (12 degree of freedom) by reducing the model in order to obtain symmetric Cauchy force-stresses.  In fact, beginning from mid 1950, Kr\"oner tried to link the theory of static
dialocations to the Cosserat model with asymmetric force stresses. However, since
1964 it was clear to Kr\"oner \cite{Kroner,Kroener64} that the force stress $\sigma$ in such a theory is
always symmetric. The relaxed micromorphic model \cite{NeffGhibaMicroModel} reconciles Kr\"{o}ner's rejection of antisymmetric force stresses in dislocation theory with the dislocation model of Eringen and Claus \cite{Eringen_Claus69,EringenClaus,Eringen_Claus71} and it is able to fully describe rotations of the microstructure and to fit a huge class of mechanical behaviors of materials with microstructure. As far as purely mechanical models are considered in the framework of linear elasticity, the need of introducing asymmetric stresses becomes rarer, see the dicussions in \cite{NeffGhibaMicroModel}. The model of Eringen and Claus \cite{Eringen_Claus69,EringenClaus,Eringen_Claus71}  contains the linear Cosserat model \cite{Neff_ZAMM05,Neff_Jeong_Conformal_ZAMM08,Neff_JeongMMS08,Neff_Paris_Maugin09} with asymmetric force stresses upon suitable restriction.

The size effects involved in a natural way in the micromorphic models  have recently received much attention in conjuction with nano-devices and foam-like structures.  {Also
other microstructured materials, as granular assemblies are
considered to be good candidates for the exploitation of
continuum micromorphic theories. Indeed, even if in the
literature the averaged models for granular assemblies are
often looked for in the framework of classical Cauchy theory
(see e.g. \cite{Chang1,Misra1,Misra2}, it becomes clearer
that generalized continuum models are necessary to correctly
describe the mechanical behavior of such physical systems
(see e.g.\cite{YangMisra,Merkel2}).}  A geometrically nonlinear generalized continuum of micromorphic type in the sense of Eringen for the phenomenological description of metallic foams is given by Neff and Forest \cite{Neff_Forest_jel05}. Moreover, in  \cite{Neff_Forest_jel05} the authors proved the existence of  minimizers and they  identified the relevant effective material parameters.  {The modelling of growth phenomena is also a major challenge to mechanical and mathematical modeling. The question of growth in continuum growth models is examined from a rigorous mathematical approach in \cite{DiCarlo1}.}

The mathematical analysis of general micromorphic solids is well-established for infinitesimal, linear elastic models, see, for example \cite{Soos69,Hlavacek69,IesanNappa2001,Iesan2002}. The only known
existence results for the static geometrically nonlinear formulation are due to Neff
\cite{Neff_micromorphic_rse_05} and to Mariano and Modica \cite{Mariano08a}.   {In fact, Mariano and Modica [94] treat general microstructures
described by manifold-valued variables, even if they discuss
essentially what is called by Neff  in  \cite{Neff_micromorphic_rse_05}  macro-stability (two
other cases are treated  in  \cite{Neff_micromorphic_rse_05}, one leads to fractures - a
situation excluded in \cite{Mariano08a} - the other is left open). When the energy
analyzed by Mariano and Modica is reduced to micromorphic materials in
the splitted version considered by Neff \cite{Mariano08a}, their coercivity assumptions
result in more stringent than Neff's ones (the blow up of the determinant
of $\textrm{det}\, F$ a part), so they restrict the material response. However, the
direct comparison of the two existence results is not completely
straightforward.}  As for the numerical implementation, see  \cite{Mariano05} and the development in \cite{Klawonn_Neff_Rheinbach_Vanis09}. In \cite{Klawonn_Neff_Rheinbach_Vanis09}
the original problem  is decoupled into two separate problems. Corresponding domain-decomposition techniques for the subproblem related to balance of forces are investigated in \cite{Klawonn_Neff_Rheinbach_Vanis09}. On the other hand, in the classical theory of Mindlin-Eringen micromorphic elasticity, existence and uniqueness results were already established by S\'{o}os \cite{Soos69}, by Hlav{\'a}{\v c}ek \cite{Hlavacek69},  by Ie\c san and Nappa \cite{IesanNappa2001} and by Ie\c san \cite{Iesan2002} assuming that the free energy is  a pointwise positive definite quadratic form. Ie\c san \cite{Iesan2002} also gave a uniqueness result for the dynamic problem without assuming that the free energy is  a positive definite quadratic form.
Moreover, in order to study the existence of solution of the resulting  system, Hlav{\'a}{\v c}ek \cite{Hlavacek69}, Ie\c san and Nappa \cite{IesanNappa2001} and Ie\c san \cite{Iesan2002} considered the  strong anchoring boundary condition.  In contrast with the models considered until now, our
free energy of the relaxed model is {\it not uniformly pointwise positive definite} in the
control of the constitutive variables. To be more precise, let us recall that the elastic free energy from the Mindlin-Eringen micromorphic elasticity model can be written as  {(see the Sections \ref{sect-notation} and \ref{sect-mod}  for notation and for the physical significations of the quantities)}
\begin{align}
2&\widehat{\mathcal{E}}(e,\varepsilon_p,\gamma)=\langle \widehat{\mathbb{C}}.\,(\nabla u-P),(\nabla u-P)\rangle
+\langle {\mathbb{H}}. \, \sym P,\sym P\rangle+
\langle \widehat{\mathbb{L}}.\,\nabla P,\nabla P\rangle\notag\\\notag\ \ \ \ &\quad\quad\quad\quad\quad\quad+
2\langle \widehat{\mathbb{E}}.\,  \sym P ,(\nabla u-P)\rangle+
2\langle \widehat{\mathbb{F}}.\, \nabla P,(\nabla u-P)\rangle+
2\langle \widehat{\mathbb{G}}.\, \nabla P, \sym P \rangle\, ,\notag
\end{align}
where $u$ is the displacement and $P$ is the
micro-distortion,  { $\langle\cdot,\cdot\rangle$ is the standard Euclidean scalar product on $\mathbb{R}^{3\times 3}$}, the constitutive coefficients are such that
\begin{align}\widehat{\mathbb{C}}:\mathbb{R}^{3\times 3}\rightarrow \mathbb{R}^{3\times 3},\quad\quad \quad \mathbb{H}:\Sym(3)\rightarrow\Sym(3),\quad\quad \quad \widehat{\mathbb{E}}, \widehat{\mathbb{G}}: \mathbb{R}^{3\times 3} \rightarrow\Sym(3), \quad\quad \quad\widehat{\mathbb{L}}:\mathbb{R}^{3\times 3\times 3}\rightarrow \mathbb{R}^{3\times 3\times 3}\notag\,,
\end{align}
and the constitutive variables are
\begin{align}
e:=\nabla u-P, \quad \quad \quad \varepsilon_p:=\sym P, \quad \quad \quad \gamma:=\nabla P \notag.
\end{align}
The elastic free energy of our relaxed model is given by
\begin{align}\label{energyourrel}
2\,&\mathcal{E}(\varepsilon_e,\varepsilon_p,\alpha)=\underbrace{\langle \C.\, \sym(\nabla u-P), \sym(\nabla u-P)\rangle}_{\text{elastic energy}}
+ \underbrace{\langle \H.\,\sym P, \sym P\rangle}_{\text{microstrain self-energy}}+ \underbrace{\langle \L.\, \Curl P, \Curl P\rangle}_{\text{dislocation energy}},\notag
\end{align}
where
\begin{align}
&\mathbb{C}:\Sym(3)\rightarrow \Sym(3),\quad\quad \quad \mathbb{H}:\Sym(3)\rightarrow\Sym(3),\quad\quad \quad
\mathbb{L}_c:\mathbb{R}^{3\times 3}\rightarrow \mathbb{R}^{3\times 3},\notag
\end{align}
and the new set of constitutive variables is
\begin{align}
\varepsilon_e=\sym(\nabla u-P), \quad \quad \quad \varepsilon_p=\sym P, \quad \quad \quad \alpha=-\Curl P.\notag
\end{align}
The comparison of the relaxed model with the classical Mindlin-Eringen \cite{Eringen99} free energy is then achieved through observing that
\begin{align}
&\langle \widehat{\mathbb{C}}.X,X\rangle_{\mathbb{R}^{3\times 3}}:=\langle \mathbb{C}.\sym X,\sym X\rangle_{\mathbb{R}^{3\times 3}},\notag\\\notag
&\langle \widehat{\mathbb{L}}.\nabla P,\nabla P\rangle_{\mathbb{R}^{3\times 3\times 3}}:=\langle \mathbb{L}_c.\Curl P,\Curl P\rangle_{\mathbb{R}^{3\times 3}}
\end{align}
define only {\it positive semi-definite tensors $\widehat{\mathbb{C}}$ and $\widehat{\mathbb{L}}$} in terms of {\it positive definite tensors $\mathbb{C}$ and $\mathbb{L}_c$} acting on linear subspaces of $\gl(3)\cong \mathbb{R}^{3\times3}$.

We prove that the new micromorphic relaxed model  \cite{NeffGhibaMicroModel} is still well-posed, i.e. we study the continuous dependence of solution with respect to the initial data and supply terms and existence and uniqueness of the solution. These results were announced previous by \cite{NeffGhibaMicroModel}. All the results are obtained for
 a standard set of tangential boundary conditions for the micro-distortion, i.e. $P.\,\tau=0$ ($P\times\,n=0$) on $\partial \Omega$ and not the usual strong anchoring condition $P=0$ on $\partial \Omega$. The solution space for the elastic distortion and micro-distortion is only ${\rm H}(\Curl;\Omega)$ and for the macroscopic displacement $u\in {\rm H}^1(\Omega)$. For non-smooth external data we expect slip lines. Using a
fundamental identity which characterizes the conservation of the total energy associated to the solution of the dynamical problem of the relaxed micromorphic model we prove the uniqueness and the continuous dependence of the solution with respect to the initial data. These results show that the considered model is in concordance with physical reality. Then, we transform
the initial boundary value problem in an abstract evolution equation in an
appropriate Hilbert space and  we use the results of the
semigroups theory of linear operators \cite{Pazy,Vrabie} in order to
obtain the existence results. The main point in establishing the desired estimates is represented by the new coercive inequalities recently proved by Neff, Pauly and  Witsch \cite{NeffPaulyWitsch,NPW2,NPW3} and by Bauer, Neff, Pauly and Starke \cite{BNPS1,BNPS2,BNPS3} (see also \cite{LNPzamp2013}).  The results established in our paper can be easily extended to theories which include electromagnetic and thermal interactions \cite{GGI11,GalesEJMA12,Grekova05,Maugin13}.

In  \cite{MadeoNeffGhibaW} we  investigate the salient  features of the new relaxed model with respect to wave-propagation phenomena  compared with the classical Mindlin-Eringen micromorphic model \cite{Mindlin64,Eringen64,Eringen99}.  {In particular, we show that the
considered relaxed model is able to account for the
description of frequency band-gaps which are observed in
particular microstructured materials as phononic crystals and
lattice structures. In particular, such materials can inhibit
wave propagation in particular frequency ranges (band-gaps)
and could be used as an alternative to piezoelectric
materials which are used today for vibration control and
which are for this reason extensively studied in the
literature (see e.g.
\cite{Piezo1,isola2,isola1,isola6,isola3,Piezo}).} Moreover, in a forthcoming paper we will deal with the static model and consider the elliptic regularity question. The numerical treatment of our new model needs FEM-discretisations in ${\rm H}({\rm curl};\Omega)$. This will be left for future work.

\section{Notation}\label{sect-notation}

For $a,b\in\R^3$ we let $\Mprod{a}{b}_{\R^3}$  denote the scalar product on $\R^3$ with
associated vector norm $\norm{a}_{\R^3}^2=\Mprod{a}{a}_{\R^3}$.
We denote by $\R^{3\times 3}$ the set of real $3\times 3$ second order tensors, written with
capital letters.
The standard Euclidean scalar product on $\R^{3\times 3}$ is given by
$\Mprod{X}{Y}_{\R^{3\times3}}=\tr({X Y^T})$, and thus the Frobenius tensor norm is
$\norm{X}^2=\Mprod{X}{X}_{\R^{3\times3}}$. In the following we omit the index
$\R^3,\R^{3\times3}$. The identity tensor on $\R^{3\times3}$ will be denoted by $\id$, so that
$\tr({X})=\Mprod{X}{\id}$.
We let $\Sym$  denote the set of symmetric tensors. We adopt the usual abbreviations of Lie-algebra theory, i.e.,
 $\so(3):=\{X\in\mathbb{R}^{3\times3}\;|X^T=-X\}$ is the Lie-algebra of  skew symmetric tensors
and $\sL(3):=\{X\in\mathbb{R}^{3\times3}\;|\, \tr({X})=0\}$ is the Lie-algebra of traceless tensors.
 For all $X\in\mathbb{R}^{3\times3}$ we set $\sym X=\frac{1}{2}(X^T+X)\in\Sym$, $\skew X=\frac{1}{2}(X-X^T)\in \so(3)$ and the deviatoric part $\dev X=X-\frac{1}{3}\;\tr({X})\,\id\in \sL(3)$  and we have
the orthogonal Cartan-decomposition  of the Lie-algebra $\gl(3)$
\begin{align}
\gl(3)&=\{\sL(3)\cap \Sym(3)\}\oplus\so(3) \oplus\mathbb{R}\!\cdot\! \id,\notag\\
X&=\dev \sym X+ \skew X+\frac{1}{3}\tr(X)\!\cdot\! \id\,.
\end{align}
By $\Co$ we denote infinitely
differentiable functions with compact support in $\Omega$. We employ the standard notation of Sobolev spaces, i.e.
$\Lp{2},\Hmp{1}{2},H_0^{1,2}(\Omega)$, which we use indifferently for
scalar-valued functions
as well as for vector-valued and tensor-valued functions. Throughout this paper (when we do not specify else) Latin subscripts take the values $1,2,3$. Typical conventions for differential
operations are implied such as comma followed
by a subscript to denote the partial derivative with respect to
 the corresponding cartesian coordinate, while $t$ after a comma denotes the partial derivative with respect to the time.
 The usual Lebesgue spaces of square integrable functions, vector or tensor fields on $\Omega$ with values in $\mathbb{R}$, $\mathbb{R}^3$ or $\mathbb{R}^{3\times 3}$, respectively will be denoted by $L^2(\Omega)$. Moreover, we introduce the standard Sobolev spaces {\cite{Adams75,Raviart79,Leis86}}
\begin{align}
&{\rm H}^1(\Omega)=\{u\in L^2(\Omega)\, |\, {\rm grad}\, u\in L^2(\Omega)\}, \ \ \ {\rm grad}=\nabla\, ,\notag\\\notag
&\ \ \ \ \ \ \ \ \ \ \ \ \ \ \ \|u\|^2_{{\rm H}^1(\Omega)}:=\|u\|^2_{L^2(\Omega)}+\|{\rm grad}\, u\|^2_{L^2(\Omega)}\, ,\\
&{\rm H}({\rm curl};\Omega)=\{v\in L^2(\Omega)\, |\, {\rm curl}\, v\in L^2(\Omega)\}, \ \ \ {\rm curl}=\nabla\times\, ,\\\notag
&\ \ \ \ \ \ \ \ \ \ \ \ \ \ \ \ \ \ \|v\|^2_{{\rm H}({\rm curl};\Omega)}:=\|v\|^2_{L^2(\Omega)}+\|{\rm curl}\, v\|^2_{L^2(\Omega)}\, ,
\end{align}
of functions $u$ or vector fields $v$, respectively.

Furthermore, we introduce their closed subspaces $H_0^1(\Omega)$, and ${\rm H}_0({\rm curl};\Omega)$ as completion under the respective graph norms of the scalar valued space $C_0^\infty(\Omega)$, the set of smooth functions with compact support in $\Omega$. Roughly speaking, $H_0^1(\Omega)$ is the subspace of functions $u\in H^1(\Omega)$ which are
zero on $\partial \Omega$, while ${\rm H}_0({\rm curl};\Omega)$ is the subspace of vectors $v\in{\rm H}({\rm curl};\Omega)$ which are normal at $\partial \Omega$  {(see \cite{Raviart79,Leis86,NeffPaulyWitsch,NPW2,NPW3})}. For vector fields $v$ with components in ${\rm H}^{1}(\Omega)$ and tensor fields $P$ with rows in ${\rm H}({\rm curl}\,; \Omega)$, i.e.,
\begin{align}
v=\left(
  \begin{array}{c}
    v_1 \\
    v_2 \\
    v_3 \\
  \end{array}
\right)\, , v_i\in {\rm H}^{1}(\Omega),
\ \quad
P=\left(
  \begin{array}{c}
    P_1^T \\
    P_2^T \\
    P_3^T \\
  \end{array}
\right)\, \quad P_i\in {\rm H}({\rm curl}\,; \Omega)\,
\end{align}
we define
\begin{align}
 {\rm Grad}\,v=\left(
  \begin{array}{c}
   {\rm grad}^T\,  v_1 \\
    {\rm grad}^T\, v_2 \\
    {\rm grad}^T\, v_3 \\
  \end{array}
\right)\, ,
\ \ \ \ {\rm Curl}\,P=\left(
  \begin{array}{c}
   {\rm curl}^T\, P_1 \\
    {\rm curl}^T\,P_2 \\
    {\rm curl}^T\,P_3 \\
  \end{array}
\right)\, , {
\ \ \ \ {\rm Grad}\,P=\left( {\rm Grad}\, P_1,{\rm Grad}\, P_2,{\rm Grad}\, P_3\right)\,.}
\end{align}

We note that $v$ is a vector field, whereas $P$, ${\rm Curl}\, P$ and ${\rm Grad}\, v$ are second order tensor fields. The corresponding Sobolev spaces will be denoted by
\begin{align}
{\rm H}({\rm Grad}\,; \Omega) \ \ \ \text{and}\ \ \ \ {\rm H}({\rm Curl}\,; \Omega)\, .
\end{align}

We recall that if $\mathbb{C}$ is a  fourth order tensor and $X\in \mathbb{R}^{3\times 3}$, then $\C. X\in \mathbb{R}^{3\times 3}$ with the components
\begin{align}
(\C.X)_{ij}=\sum\limits_{k=1}^{3}\sum\limits_{l=1}^{3}\mathbb{C}_{ijkl}X_{kl}\, .
\end{align}

\section{Formulation of the relaxed micromorphic continuum  model}\label{sect-mod}\setcounter{equation}{0}

We consider a micromorphic continuum which occupies a bounded domain $\Omega$ and is bounded by the piecewise smooth
surface $\partial \Omega$. Let $T>0$ be a given time. The motion of the body is referred to a fixed system of rectangular Cartesian axes $Ox_i$, $(i=1,2,3)$.
The micro-distortion (plastic distortion) $P=(P_{ij}):\Omega\times [0,T]\rightarrow \mathbb{R}^{3 \times 3}$  describes the substructure of the material which can rotate, stretch, shear and shrink, while \linebreak $u=(u_i) :\Omega\times [0,T]\rightarrow  \mathbb{R}^3$  is the displacement of the macroscopic material points.

The quantities involved in our new relaxed micromorphic continuum  model have the following physical signification:
\begin{itemize}
\item $(u,P)$ are the {\it kinematical variables},
\item $u$ is the displacement vector (translational degrees of freedom),
\item $P$ is the micro-distortion tensor (plastic distortion, second order, non-symmetric),
\item ${\sigma}$ is   the  Cauchy stresses (second order, symmetric),
\item $s$ is the microstress tensor (second order, symmetric),
\item ${m}$ is the moment stress  tensor (micro-hyperstress tensor, third order,  in general non-symmetric),
\item $f$ is the body force,
\item $M$ is the body moment tensor (second order, non-symmetric),
\item $e:=\nabla u-P$ is the elastic distortion (relative distortion, second order, non-symmetric),
\item $\varepsilon_e:=\sym e=\sym(\nabla u-P)$ is the elastic strain tensor (second order, symmetric),
\item $\varepsilon_p:=\sym P $ is the micro-strain tensor (plastic strain, second order, symmetric),
\item  $\alpha:=\Curl e=-\Curl P$ the micro-dislocation tensor (second order).
\end{itemize}

We consider here a relaxed version of the classical micromorphic model with {\it} $\sigma$ symmetric  and drastically reduced numbers of constitutive coefficients. More precisely, our model is a subset of the  { classical  micromorphic  model in which we allow the usual micromorphic  tensors} \cite{Eringen99} to become positive-semidefinite only \cite{NeffGhibaMicroModel}. The proof of the well-posedness of this model necessitates the {\it application of new mathematical tools}   \cite{NeffPaulyWitsch,NPW2,NPW3,BNPS1,BNPS2,BNPS3,LNPzamp2013}. The curvature dependence is reduced to a dependence only on the micro-dislocation tensor $\alpha:=\Curl e=-\Curl P\in\mathbb{R}^{3\times 3}$ instead of
$\gamma=\nabla P\in\mathbb{R}^{27}=\mathbb{R}^{3\times 3\times 3}$ and the local response is reduced to a dependence on the symmetric part of the elastic distortion (relative distortion) $\varepsilon_e=\sym e=\sym(\nabla u-P)$, while the {\it full kinematical degrees of freedom} for $u$ and $P$ are kept, notably {\it rotation of the microstructure remains possible}.

Our new set of {\it independent constitutive variables} for the relaxed micromorphic model is thus
\begin{align}
\varepsilon_e=\sym(\nabla u-P), \quad \quad \quad \varepsilon_p=\sym P, \quad \quad \quad \alpha=-\Curl P.
\end{align}

The system  of partial differential equations  which corresponds to this special linear anisotropic micromorphic continuum is derived from the following free energy
\begin{align}\label{energyourrel}
\quad 2\,\mathcal{E}(\varepsilon_e,\varepsilon_p,\alpha)&=\langle \C.\, \varepsilon_e, \varepsilon_e\rangle
+ \langle \H.\,\varepsilon_p, \varepsilon_p\rangle+ \langle \L.\, \alpha, \alpha\rangle\\
&=\underbrace{\langle \C.\, \sym(\nabla u-P), \sym(\nabla u-P)\rangle}_{\text{elastic energy}}
+ \underbrace{\langle \H.\,\sym P, \sym P\rangle}_{\text{microstrain self-energy}}+ \underbrace{\langle \L.\, \Curl P, \Curl P\rangle}_{\text{dislocation energy}},\notag
\\\notag
&\hspace{-1.7cm}\sigma=D_{\varepsilon_e} \, \mathcal{E}(\varepsilon_e,\varepsilon_p,\alpha)\in \Sym(3),\quad \quad
s=D_{\varepsilon_p} \, \mathcal{E}(\varepsilon_e,\varepsilon_p,\alpha),\in \Sym(3),\quad \quad
m=D_{\alpha} \, \mathcal{E}(\varepsilon_e,\varepsilon_p,\alpha)\in \mathbb{R}^{3\times3},
\end{align}
where  $\C\!:\!\Omega\rightarrow L(\mathbb{R}^{3 \times 3},\mathbb{R}^{3 \times 3})$, $\L\!:\!\Omega\rightarrow L(\mathbb{R}^{3 \times 3},\mathbb{R}^{3 \times 3})$ and $\H\!:\!\Omega\rightarrow L(\mathbb{R}^{3 \times 3},\mathbb{R}^{3 \times 3})$
are fourth order elasticity tensors, positive definite and functions of class $C^1(\Omega)$.

For the rest of the paper we assume  that the constitutive coefficients
\begin{align}
&\mathbb{C}:\Sym(3)\rightarrow \Sym(3),\quad\quad \quad \mathbb{H}:\Sym(3)\rightarrow \Sym(3),\quad\quad \quad
\mathbb{L}_c:\mathbb{R}^{3\times 3}\rightarrow \mathbb{R}^{3\times 3}
\end{align}
have the following symmetries
\begin{align}\label{simetries}
&\C_{ijrs}=\C_{rsij}=\C_{jirs},\quad\quad\quad \quad \H_{ijrs}=\H_{rsij}=\H_{jirs}, \quad\quad\quad \quad {(\L)}_{ijrs}={(\L)}_{rsij}\,.
\end{align}

We assume that the  fourth order elasticity tensors $\C$, $\L$ and $\H$ are positive definite.
Then, there are  positive numbers ${c_M}$, ${c_m}$ (the maximum and minimum elastic moduli for $\C$), ${(L_c)}_M$, ${(L_c)}_m$ (the maximum and minimum moduli for $\L$) and $h_M$, $h_m$ (the maximum and minimum moduli for $\H$) such that
\begin{align}\label{posdef}
{c_m}\|X \|^2\leq \langle\,\C. X,X\rangle\leq {c_M}\|X \|^2\,\quad  \ \ &\text{for all }\ \ X\in\Sym(3),\notag
\\
{(L_c)}_m\|X \|^2\leq \langle \L. X,X\rangle\leq {(L_c)}_M\| X\|^2\,\quad \ \ &\text{for all }\ \ X\in\mathbb{R}^{3\times 3},
\\
h_m\|X \|^2\leq \langle \H. X,X\rangle\leq h_M\| X\|^2\, \quad \ \ &\text{for all }\ \ X\in\Sym(3).\notag
\end{align}
Further we  assume, without loss of generality, that ${c_M}$, ${c_m}$, ${(L_c)}_M$, $h_M$, $h_m$ and ${(L_c)}_m$ are constants.

We introduce the action functional of the considered system to be defined as
\begin{align}
\Lambda=\int_0^T\int_\Omega(\mathcal{K}+\varpi-\mathcal{E})\,dv\,dt,\quad\quad \mathcal{K}=\frac{1}{2}\|u_{,t}\|^2+ \frac{1}{2}\|P_{,t}\|^2\,\quad\quad  \varpi=\langle f,{u}_{,t}\rangle +\langle M,{P}_{,t}\rangle\, ,
\end{align}
where $\mathcal{K}$ and $\varpi$ are the kinetic energy and the work done by external loads on the body, \linebreak respectively,  $f :\Omega\times [0,T]\rightarrow  \mathbb{R}^3$ describes the body force and $M:\Omega\times [0,T]\rightarrow \mathbb{R}^{3 \times 3}$ describes the external body moment.

We consider the weaker boundary conditions
\begin{align} \label{bc}
{u}({x},t)=0, \ \ \  \text{and the {\it tangential condition}}  \quad {P}_i({x},t)\times n(x) =0, \ \ \ i=1,2,3, \ \ \
\ ({x},t)\in\partial \Omega\times [0,T],
\end{align}%
where $\times$ denotes the vector product, $n$ is the unit outward normal vector at the surface $\partial \Omega$, $P_i, i = 1, 2, 3$ are
the rows of $P$. The model is driven by nonzero initial conditions
\begin{align}\label{ic}
&{u}({x},0)={u}_0(
x),\quad \quad \quad\dot{u}({x},0)=\dot{u}_0(
x),\quad \quad \quad
{P}({x},0)={P}_0(
x), \quad \quad \quad\dot{P}({x},0)=\dot{P}_0(
x),\ \ \text{\ \ }{x}\in \overline{\Omega},
\end{align}%
where  ${u}_0, \dot{u}_0, {P}_0$ and $\dot{P}_0$ are prescribed functions.

\begin{remark}Since $P$ is determined in  $ {\rm H}({\rm Curl}\,; \Omega)$, in our relaxed model
the only possible description of boundary value is in terms of tangential traces, i.e.
$P.\tau=0$ for all tangential vectors $\tau$ at $\partial \Omega$. This follows from the standard theory of the  ${\rm H}({\rm Curl}\,; \Omega)$-space.
\end{remark}

Imposing the first variation of the action functional to be zero (Hamilton-Kirchhoff principle), integrating by parts a suitable
number of times and considering arbitrary variations $\delta u$ and $\delta P$ of the basic kinematical fields, we obtain that the
system of partial differential equations  of our relaxed  micromorphic continuum  model is
\begin{align}\label{eqrelax}
u_{,tt}&=\dvg[ \C. \sym(\nabla u-P)]+f\, , \quad\quad \quad\quad\quad\quad\quad\quad\quad\quad\quad\quad\quad \ \ \ \ \  \text{balance of forces}, \\\notag
P_{,tt}&=- \crl[ \L.\crl\, P]+\C. \sym (\nabla u-P)-\H. \sym P+M, \quad\quad \ \  \text{balance of moment stresses}\, \ \ \
\end{align}
in $  \Omega\times [0,T]$. For simplicity, the system \eqref{eqrelax} is considered in a normalized form.

Our new approach, in marked contrast to classical asymmetric micromorphic models, features
a {\it symmetric Cauchy stress tensor} $\sigma=\C. \sym(\nabla u-P)$. Therefore, the linear Cosserat approach (\cite{Neff_ZAMM05}: $\mu_c>0$) is excluded here.

The relaxed formulation considered  in the present  paper still shows size effects and smaller samples are relatively stiffer.  It is clear to us that for this reduced model of relaxed micromorphic elasticity {\it unphysical  effects of singular stiffening behaviour for small sample sizes} (``bounded stiffness", see \cite{Neff_JeongMMS08}) {\it cannot appear}. In case of the isotropic Cosserat model this is only true for a reduced curvature energy depending only on $\|\dev \sym \Curl P\|^2$, see the discussion in \cite{Neff_JeongMMS08,Neff_Paris_Maugin09}.

 In contrast with  the 7+11 parameters isotropic Mindlin and Eringen model \cite{Mindlin64,Eringen64,EringenSuhubi2}, we have altogether only  seven parameters  $\mu_e,\lambda_e,  \mu_h, \lambda_h,\alpha_1, \alpha_2$, $\alpha_3$.  For isotropic materials, our system \label{eq} reads
\begin{align}\label{eqiso}
 u_{,tt}&=\dvg\,\sigma+f\, ,\\\notag
P_{,tt}&=-\Curl m+\sigma-s+M\,  \ \ \ \text {in}\ \ \  \Omega\times [0,T].
\end{align}
where
\begin{align}
\sigma&= 2\mu_e \sym(\nabla u-P)+\lambda_e \tr(\nabla u-P){\cdp} \id,\notag\\
m&=\alpha_1 \dev\sym \Curl P+\alpha_2 \skew \Curl P +\alpha_3\, \tr(\Curl P){\cdp} \id,\\
s&=2\mu_h \sym P+\lambda_h \tr (P){\cdp} \id\,. \notag
\end{align}
Thus, we obtain the complete system of linear partial differential equations in terms of the kinematical unknowns $u$ and $P$
\begin{align}\label{eqisoup}
 u_{,tt}&=\dvg[2\mu_e \sym(\nabla u-P)+\lambda_e \tr(\nabla u-P){\cdp} \id]+f\, ,\\\notag
P_{,tt}&=-\Curl [\alpha_1 \dev\sym \Curl P+\alpha_2 \skew \Curl P +\alpha_3\, \tr(\Curl P){\cdp} \id]\\&\quad\ +2\mu_e \sym(\nabla u-P)+\lambda_e \tr(\nabla u-P){\cdp} \id-2\mu_h \sym P-\lambda_h \tr (P){\cdp} \id+M\,  \ \ \ \text {in}\ \ \  \Omega\times [0,T].\notag
\end{align}

In this model, {\it the asymmetric parts of $P$} are entirely due only to {\it moment stresses} and {\it applied body moments}\,! In this sense, the {\it macroscopic} and {\it microscopic scales} are neatly {\it separated}.

The positive definiteness required for  the tensors $\mathbb{C}$, $\mathbb{H}$ and $\mathbb{L}_c$ implies for isotropic materials the following restriction upon the parameters $\mu_e,\lambda_e, \mu_h, \lambda_h, \alpha_1, \alpha_2$ and $\alpha_3$
\begin{align}\label{condpara}
\mu_e>0,\quad\quad  2\mu_e+3\lambda_e>0, \quad\quad  \mu_h>0, \quad\quad  2\mu_h+3\lambda_h>0, \quad\quad  \alpha_1>0,
\quad\quad  \alpha_2>0, \quad\quad  \alpha_3>0.
\end{align}
Therefore, positive definiteness for our isotropic model does not involve extra nonlinear side conditions \cite{Eringen99,Smith}.

If, by abuse of our notation by neglect of our guiding assumption, we add the anti-symmetric term \linebreak $2\mu_c\skew(\nabla u-P)$ in the expression of the Cauchy stress tensor $\sigma$, where $\mu_c\geq0$ is the {\it Cosserat couple modulus},  then our analysis works for $\mu_c\geq0$. The model in which $\mu_c>0$ is the isotropic Eringen-Claus  model for dislocation dynamics \cite{Eringen_Claus69,EringenClaus,Eringen_Claus71} and it is derived from the following free energy
\begin{align}\label{XXXX}
\mathcal{E}(e,\varepsilon_p,\alpha)&=\mu_e \|\sym (\nabla u-P)\|^2+\mu_c\|\skew(\nabla u-P)\|^2+ \frac{\lambda_e}{2}\, [\tr(\nabla u-P)]^2+\mu_h \|\sym  P\|^2+ \frac{\lambda_h}{2} [\tr\,(P)]^2\notag\\&
 \quad \quad  +\frac{\alpha_1}{2}\| \dev\sym \Curl P\|^2 +\frac{\alpha_2}{2}\| \skew \Curl P\|+ \frac{\alpha_3}{2}\, \tr(\Curl P)^2.
 \end{align}
  For $\mu_c>0$ and if the other inequalities \eqref{condpara} are satisfied, the existence and uniqueness follow along the classical lines. There is no need for any new integral inequalities. To the sake of simplicity, we only present in the present paper well-posedness results for the relaxed model. These results still hold for the complete model and can be generalized with some additional calculations.

 For the mathematical treatment of the linear relaxed model there arises the need for new integral type inequalities which we present in the next section. Using the new results established by Neff, Pauly and  Witsch \cite{NeffPaulyWitsch,NPW2,NPW3} and by Bauer, Neff, Pauly and Starke \cite{BNPS1,BNPS2,BNPS3} we are now able to manage also  energies depending on the dislocation energy and  having symmetric Cauchy stresses.

\section{New Poincar\'{e} and Korn type estimates}\setcounter{equation}{0}

In potential theory use is made of Poincar\'{e}'s inequality, that is
\begin{align}\label{Poincare}
\|u\|_{L^2(\Omega)}\leq c_p \|{\rm grad}\, u\|_{L^2(\Omega)}\, ,
\end{align}
for all functions $u\in H_0^1(\Omega)$ with some constants $c_p>0$, to bound a scalar potential in terms of its gradient.

In linearized elasticity theory Korn's inequality is used, that is
\begin{align}\label{Korn}
\|{\rm grad}\,u\|_{L^2(\Omega)}\leq c_k \|\sym {\rm grad}\, u\|_{L^2(\Omega)}\, ,
\end{align}
for all functions $u\in H_0^1(\Omega)$ with some constants $c_k>0$, for bounding the deformation of an elastic medium in terms of the symmetric strains.

In electro-magnetic theory the Maxwell inequality, that is
\begin{align}\label{Maxwell}
\|u\|_{L^2(\Omega)}\leq c_m( \|{\rm curl}\, u\|_{L^2(\Omega)}+\|{\rm div}\, u\|_{L^2(\Omega)})\, ,
\end{align}
for all functions $u\in H_0({\rm curl}; \Omega)\cap H({\rm div};\Omega)$ with some constants $c_m>0$, is used to bound the electric and magnetic field in terms of the electric charge and current density, respectively.

In \cite{NeffPaulyWitsch,NPW2,NPW3}, for tensor fields $P\in {\rm H}({\rm Curl}\, ; \Omega)$ the following seminorm $|\|\cdot|\|$ is defined
\begin{align}
\||P\||^2:=\|\sym P\|^2_{L^2(\Omega)}+\|{\rm Curl}\, P\|^2_{L^2(\Omega)}\, .
\end{align}

From \cite{NeffPaulyWitsch,NPW2,NPW3} we have the following result:
\begin{theorem}
There exists a constant $\hat{c}$ such that
\begin{align}\label{Neff1}
\| P\|_{L^2(\Omega)}\leq \hat{c}\||P\||\, ,
\end{align}
for all $P\in {\rm H}({\rm Curl}\, ; \Omega)$ with vanishing restricted tangential trace on $\partial \Omega$, i.e. $P.\tau=0$ on $\partial \Omega$ .
\end{theorem}

Moreover, we have
\begin{theorem}\label{wdn}
On ${\rm H}_0({\rm Curl}\, ; \Omega)$  the norms $\|\cdot\|_{{\rm H}({\rm Curl}\, ; \Omega)}$ and $|\|\cdot|\|$ are equivalent. In particular,
$|\|\cdot|\|$ is a norm on ${\rm H}_0({\rm Curl}\, ; \Omega)$ and there exists a positive constant $c$, such that
\begin{align}\label{Neff2}
c\,\|P\|_{{\rm H}({\rm Curl}\, ; \Omega)}\leq \||P\||\, ,
\end{align}
for all $P\in {\rm H}_0({\rm Curl}\, ; \Omega)$.
\end{theorem}

Moreover,  in a forthcoming paper \cite{BNPS2} (see also \cite{BNPS1,NeffPaulyWitsch}) the following results are proved:

\begin{theorem}\label{BNPS2}
There exists a positive constant $C_{DD}$, only depending on $\Omega$, such that for all $P\in{\rm H}_0({\rm Curl}\, ; \Omega)$ the following estimate hold:
\begin{align}\label{BNPSe2}
\| \Curl P\|_{L^2(\Omega)}\leq C_{DD}\,\|\dev  \Curl P\|_{L^2(\Omega)}\, .
\end{align}
\end{theorem}

\begin{theorem}\label{BNPS}
There exists a positive constant $C_{DSDC}$, only depending on $\Omega$, such that for all $P\in{\rm H}_0({\rm Curl}\, ; \Omega)$ the following estimate hold:
\begin{align}\label{BNPSe}
\| P\|_{L^2(\Omega)}\leq C_{DSDC}(\|\dev \sym P\|^2_{L^2(\Omega)}+\|\dev \Curl P\|^2_{L^2(\Omega)})\, .
\end{align}
\end{theorem}

\begin{corollary}\label{BNPS3}
For all $P\in{\rm H}_0({\rm Curl}\, ; \Omega)$ the following estimate hold:
\begin{align}\label{BNPSe3}
\| P\|_{L^2(\Omega)}+\| \Curl P\|_{L^2(\Omega)}\leq (C_{DSDC} +C_{DD})(\|\dev \sym P\|^2_{L^2(\Omega)}+\|\dev \Curl P\|^2_{L^2(\Omega)})\, .
\end{align}
\end{corollary}

We have to remark that the above corollary proves that on ${\rm H}_0({\rm Curl}\, ; \Omega)$  the norms $\|\cdot\|_{{\rm H}({\rm Curl}\, ; \Omega)}$ and $\|\|\cdot\|\|$ are equivalent.

\begin{theorem}\label{BNPS4}
There exists a positive constant $C_{DSG}$, only depending on $\Omega$, such that for all $u\in{\rm H}_0^1(\Omega)$ the following estimate hold:
\begin{align}\label{BNPSe4}
\| \nabla u\|_{L^2(\Omega)}\leq C_{DSG}\|\dev\sym \nabla u\|_{L^2(\Omega)}\, .
\end{align}
\end{theorem}

The estimates given by the above theorems will be essential in the study of our relaxed linear micromorphic model.

\section{Conservation law, uniqueness, continuous dependence and \\existence}
\setcounter{equation}{0}
\subsection{Energy conservation}
In this subsection we establish a
fundamental identity which characterize the conservation of the total energy associated to the solution of the dynamic problem $(\mathcal{P})$ defined by the equations \eqref{eqrelax}, the boundary conditions \eqref{bc} and the initial conditions \eqref{ic}.
Let us consider a solution  $\{u,P\}$ of the problem  $(\mathcal{%
P})$ corresponding to the given data
$\mathcal{I}=\{f,M, {u}_0, \dot{u}_0, {P}_0, \dot{P}_0\}.$

We define the total energy
\begin{align}
2E(t)=&\dd\int_\Omega\bigg(\|u_{,t}\|^2+ \|P_{,t}\|^2+ \langle \C.\, \sym(\nabla u-P), \sym(\nabla u-P)\rangle
\vspace{2mm}\notag\\&\ \ \ \ \ \ \quad\quad + \langle \H.\,\sym P, \sym P\rangle+ \langle \L.\, \Curl P, \Curl P\rangle\bigg)dv\, ,
\end{align}
and the power function
\begin{align}
 \Pi(t)&=\dd\int_\Omega(\langle f,{u}_{,t}\rangle +\langle M,{P}_{,t}\rangle%
)dv\, .
\end{align}

\begin{lemma}\label{lemaceg} {\rm (Conservation law)} Let $%
\{u,P\}$
be a solution of the problem $(\mathcal{P})$ corresponding to the loads \break $\mathcal{I}=\{f,M, {u}_0, \dot{u}_0, {P}_0, \dot{P}_0\}.$ Then, for every time $t\in [0,T]$, we have
\begin{equation}\label{cons}
E(t)=E(0)+\int_0^t\, \Pi(s)ds.
\end{equation}
\end{lemma}
\textit{Proof.} First of all, let us recall the identities
\begin{align}
&{\rm div} (\psi A)=\langle A, {\rm grad} \, \psi\rangle+\psi\, {\rm div}\,  A\, ,\\
&{\rm div}\,  (A\times B)=\langle B, \curl \, A\rangle-\langle A, \curl\,  B\rangle\, ,\notag
\end{align}
for all $C^1$-functions $\psi:\Omega\rightarrow \mathbb{R}$ and $A,B:\Omega\rightarrow \mathbb{R}^{3}$, where $\times$ is the cross product. Hence
\begin{align}\label{formule}
&{\rm div} (\varphi_i Q_i)=\langle Q_i, \nabla \, \varphi_i\rangle+\varphi_i\, {\rm div}\,  Q_i\,\quad  \ \ \text{not summed} ,\\
&{\rm div}\,  (R_i\times S_i)=\langle S_i, \curl \, R_i\rangle-\langle R_i, \curl\,  S_i\rangle\, \quad  \ \ \text{not summed} ,\notag
\end{align}
for all $C^1$-functions $\varphi_i:\Omega\rightarrow \mathbb{R}$ and $Q_i,P_i,S_i:\Omega\rightarrow \mathbb{R}^{3}$, where $\varphi_i$ are the components of the vector $\varphi$ and $Q_i, P_i,S_i$ are the rows of the matrix $Q$, $P$ and $S$, respectively. We choose
\begin{align}
\varphi=u_{,t},\quad\quad Q=\C. \sym(\nabla u-P)
\end{align}
and we obtain
\begin{align}
&{\rm div} (u_{i,t} [\C. \sym(\nabla u-P)]_i)=\langle [\C. \sym(\nabla u-P)]_i, \nabla \, u_{i,t}\rangle+u_{i,t}\, {\rm div}\,  [\C. \sym(\nabla u-P)]_i\,
 \ \ \text{not summed}.
\end{align}
This leads to
\begin{align}
&\sum\limits_{i=1}^3u_{i,t}\, {\rm div}\,  [\C. \sym(\nabla u-P)]_i=\sum\limits_{i=1}^3{\rm div} (u_{i,t} [\C. \sym(\nabla u-P)]_i)-\sum\limits_{i=1}^3\langle [\C.\sym(\nabla u-P)]_i, \nabla \, u_{i,t}\rangle\, .
\end{align}
Thus
\begin{align}\label{fdiv}
&\langle  {\rm Div}\,  [\C. \sym(\nabla u-P)], u_{,t}\rangle=\sum\limits_{i=1}^3{\rm div} (u_{i,t} [\C. \sym(\nabla u-P)]_i)-\langle \C. \sym(\nabla u-P), \sym \nabla u_{,t}\rangle\, .
\end{align}
If we take in \eqref{formule}
\begin{align}
R_i=[\L. \Curl P]_i, \quad\quad S_i=P_i\, ,
\end{align}
we have
\begin{align}
&\sum\limits_{i=1}^3{\rm div}\,  ([\L.\curl P]_i\times P_{i,t})=\sum\limits_{i=1}^3\langle P_{i,t}, \curl \, [\L. \Curl P]_i\rangle-\sum\limits_{i=1}^3\langle [\L. \Curl P]_i, \curl\,  P_{i,t}\rangle\, .\notag
\end{align}
Hence, we obtain
\begin{align}\label{fcurl}
&\langle P_{,t}, \Curl \, (\L. (\Curl \,  P))\rangle=\sum\limits_{i=1}^3{\rm div}\,  ([\L. \Curl P]_i\times P_{i,t})+\langle \L. \Curl\, P, \Curl\, P_{,t}\rangle\, .
\end{align}
Using \eqref{eqrelax}, \eqref{fdiv} and \eqref{fcurl} we have
\begin{align}
\langle u_{,tt}, u_{,t}\rangle+\langle P_{,tt}, P_{,t}\rangle&=\langle\dvg  (\C. \sym(\nabla u-P)), u_{,t}\rangle+\langle f, u_{,t}\rangle\, \\\notag
&-\langle\crl(\L. \crl( P)), P_{,t}\rangle+\langle \C. \sym (\nabla u-P), P_{,t}\rangle-\langle \H\, \sym P,  P_{,t}\rangle+\langle M,  P_{,t}\rangle\,\\\notag
=&\sum\limits_{i=1}^3{\rm div} (u_{i,t} [\C. \sym(\nabla u-P)]_i)-\langle \C. \sym(\nabla u-P), \sym \nabla u_{,t}\rangle+\langle f, u_{,t}\rangle\, \\\notag
&-\sum\limits_{i=1}^3{\rm div}\,  [(\L.\, \Curl P)_i\times P_{i,t}]-\langle \L. \Curl\, P, \Curl\, P_{,t}\rangle\\\notag
&+\langle \C. \sym (\nabla u-P), P_{,t}\rangle-\langle \H. \sym P,  P_{,t}\rangle+\langle M,  P_{,t}\rangle\,\notag\\\notag
=&-\langle \C. \sym(\nabla u-P), \sym(\nabla u_{,t}-P_{,t})\rangle-\langle \L. \Curl\, P, \Curl\, P_{,t}\rangle-\langle \H. \sym P,  \sym P_{,t}\rangle\, \\\notag
&+\sum\limits_{i=1}^3{\rm div} (u_{i,t} [\C. \sym(\nabla u-P)]_i)+\sum\limits_{i=1}^3{\rm div}\,  ( P_{i,t}\times [\L. \Curl P]_i)\\\notag
&+\langle f, u_{,t}\rangle+\langle M,  P_{,t}\rangle\,.\notag
\end{align}
Hence, using the symmetries \eqref{simetries}, we have
\begin{align}
\frac{1}{2}\frac{\partial }{\partial t}\bigg(&\|u_{,t}\|^2+ \|P_{,t}\|^2+ \langle \C. \sym(\nabla u-P), \sym(\nabla u-P)\rangle
\vspace{2mm}\notag\\&+ \langle \H.\sym P, \sym P\rangle+ \langle \L. \Curl P, \Curl P\rangle\, \bigg) \\\notag
&=\sum\limits_{i=1}^3{\rm div} (u_{i,t} [\C. \sym(\nabla u-P)]_i)+\sum\limits_{i=1}^3{\rm div}\,  ( P_{i,t}\times [\L. \Curl P]_i)\\\notag
&+\langle f, u_{,t}\rangle+\langle M,  P_{,t}\rangle.\notag
\end{align}
Therefore, using the divergence theorem, it follows that
\begin{align}
\frac{d}{dt}E(t)=&\dd\int_{\partial \Omega}\bigg(\sum\limits_{i=1}^3 \langle [\C.\sym(\nabla u-P)]_iu_{i,t},n\rangle+\sum\limits_{i=1}^3\langle P_{i,t}\times [\L. \Curl P]_i,n\rangle\bigg)da\\\notag
&+\dd\int_\Omega(\langle f,{u}_{,t}\rangle +\langle M,{P}_{,t}\rangle%
)dv\\\notag
=&\dd\int_{\partial \Omega}\bigg(\sum\limits_{i=1}^3 \langle [\C. \sym(\nabla u-P)]_iu_{i,t},n\rangle+\sum\limits_{i=1}^3(\langle [\L. \Curl P]_i, n\times P_{i,t} \rangle\bigg)da\\\notag
&+\dd\int_\Omega(\langle f,{u}_{,t}\rangle +\langle M,{P}_{,t}\rangle%
)dv\, ,
\end{align}%
so that, in view of  the  boundary conditions $u=0,\ P.\tau=0$ on $\partial \Omega$ and by integration over $[0,t]$, the proof is complete.\hfill $\Box$

 \subsection{Continuous dependence of solution and uniqueness}

Throughout this section we study the continuous dependence of
solution of the problem $(\mathcal{P})$ with respect to the initial  and the
body loads. To this aim, let us first prove the following lemma:
\begin{lemma}\label{lemaii} There exists a constant $a_1$ such that
\begin{align}
a_1(\| \sym \nabla u\|^2+\|\sym P\|^2)\leq \langle \C. \sym(\nabla u-P), \sym(\nabla u-P)\rangle+ \langle \H.\sym P, \sym P \rangle
\end{align}
for all $u\in H^1(\Omega)$ and $P\in H(\Curl;\Omega)$.
\end{lemma}

{\it Proof.} We start the proof with the remark that the arithmetic-geometric inequality and the positivity of $\C$ imply
\begin{align}\label{ineqmed}
\langle \C. \sym(\nabla u-P), \sym(\nabla u-P)\rangle&\geq {c_m} \|\sym(\nabla u-P)\|^2\notag
= {c_m}\langle \sym(\nabla u-P), \sym(\nabla u-P)\rangle\\\notag
&\geq {c_m}[\| \sym \nabla u \|^2+\|\sym P\|^2-2\langle \sym \nabla u,\sym P\rangle]
\\
&\geq {c_m}\bigg[(1-\delta)\| \sym \nabla u\|^2+\bigg(1 -\frac{1}{\delta}\bigg)\|\sym P\|^2\bigg]\, ,
\end{align}
for all $\delta>0$. Moreover, we have
\begin{align}
 \langle \H.\sym P, \sym P\rangle\geq h_m\,\|\sym P\|^2\,.
\end{align}
Hence, we deduce that
\begin{align}
 {c_m}(1-\delta)\| \sym \nabla u \|^2+\bigg({c_m} +h_m-\frac{{c_m}}{\delta}\bigg)\|\sym P\|^2 \leq \langle& \C. \sym(\nabla u-P), \sym(\nabla u-P)\rangle\\&\notag+ \langle \H.\sym P, \sym P \rangle.
\end{align}

If we choose $\delta$ so that
$$
\frac{{c_m} }{{c_m} +h_m}<\delta<1
$$
we have that there is a positive constant $a_1$ so that
\begin{align}
a_1\big( \| \sym \nabla u \|^2+\|\sym P \|^2\big)\leq \langle \C. \sym(\nabla u-P), \sym(\nabla u-P)\rangle+ \langle \H.\sym P , \sym P \rangle\,
\end{align}
and the proof is complete.\hfill $\Box$

To establish an estimate describing the continuous dependence upon
the
initial data we shall assume that  $\{u,P\}$ is  solution of the problem
$(\mathcal{P})$ with null boundary data and null body loads. For this type of external data system, using Lemma
 \ref{lemaceg}, we deduce the following result.

\begin{theorem}\label{cduid}{\rm (Continuous dependence upon initial data)}  Let $\{u,P\}$ %
be a solution of the problem  $(\mathcal{P})$
with the external data system $\mathcal{I}=\{0,0, {u}_0, \dot{u}_0, {P}_0, \dot{P}_0\}.$
 Then, there is a positive constant $a$ so that
\begin{equation}\label{cdi}
a\big(\|u_{,t}\|^2_{L^2(\Omega)}+ \|P_{,t}\|^2_{L^2(\Omega)}+ \| \nabla u\|^2_{L^2(\Omega)}+\|P\|^2_{L^2(\Omega)}+ \|\Curl P \|^2_{L^2(\Omega)}\big)\leq E(0), \quad \ \text{for all} \ \  t\in[0,T]\, .
\end{equation}
\end{theorem}

{\it Proof.} A direct consequence of the conservation law is
\begin{align}
\|u_{,t}\|^2+ \|P_{,t}\|^2&+ \langle \C\sym(\nabla u-P), \sym(\nabla u-P)\rangle
\vspace{2mm}\notag\\&+ \langle \H.\sym P, \sym P\rangle+ \langle \L. \Curl P, \Curl P\rangle=2\,E(0),  \quad\  \text{for all} \ \  t\in[0,T]\, .
\end{align}

Using Lemma \ref{lemaii} and the inequality
\begin{align}
\langle \L. \Curl P, \Curl P\rangle\geq {(L_c)}_m\,\|\Curl P\|^2\, ,
\end{align}
 we have
 that there is a positive constant $a_2$ so that
\begin{align}
a_2\big(\|u_{,t}\|^2_{L^2(\Omega)}+ \|P_{,t}\|^2_{L^2(\Omega)}+ \| \sym \nabla u \|^2_{L^2(\Omega)}+\|\sym P\|^2_{L^2(\Omega)}+ \|\Curl P \|^2_{L^2(\Omega)}\big)\leq E(0)\, .
\end{align}

Because  $P\in{\rm H}_0({\rm Curl}\, ; \Omega)$ and $u\in H_0^1(\Omega)$, in view of \eqref{Korn} and \eqref{Neff2}, there are the positive constants $c$ and $c_k$ \cite{NeffPaulyWitsch,NPW2,NPW3}, such that
\begin{align}
c\|P\|_{{\rm H}({\rm Curl}\, ; \Omega)}\leq \||P\||\, ,
\end{align}
and
\begin{align}
\|\nabla \,u\|_{L^2(\Omega)}\leq c_k \|\sym \nabla\, u\|_{L^2(\Omega)}\, .
\end{align}

Hence, we can find a positive constant $a$ so that
\begin{align}
a\big(\|u_{,t}\|^2_{L^2(\Omega)}+ \|P_{,t}\|^2_{L^2(\Omega)}+ \| \nabla u\|^2_{L^2(\Omega)}+\|P\|^2_{L^2(\Omega)}+ \|\Curl P \|^2_{L^2(\Omega)}\big)\leq E(0), \ \ \text{for all} \quad \  t\in[0,T]\, .\notag\text{\quad\quad\quad\quad\quad\hfill $\Box$}
\end{align}

\begin{corollary} {(\rm Uniqueness)}
Any two solutions of the problem $(\mathcal{P})$ are equal.\hfill $\Box$
\end{corollary}

Now we study the continuous data dependence of the solution upon the
supply terms  $\{f,M\}$.
\begin{theorem}\label{cdst} {\rm (Continuous dependence upon the supply terms)} %
 Let  $\{u,P\}$ be a solution of the problem $(\mathcal{P})$
corresponding to external data system
$\mathcal{I}=\{f,M, 0,0,0,0\}.$
Then, for all $t\in I$ we have
\begin{equation}
\sqrt{a}\big(\|u_{,t}\|^2_{L^2(\Omega)}+ \|P_{,t}\|^2_{L^2(\Omega)}+ \| \nabla u\|^2_{L^2(\Omega)}+\|P\|^2_{L^2(\Omega)}+ \|\Curl P \|^2_{L^2(\Omega)}\big)^{\frac{1}{2}}\leq \frac{1}{2}\int _0^t g(s)ds,
\end{equation}%
where
\begin{equation}\label{2.3.53}
g(s)=\bigg\{\dd\int_\Omega(\| f(s)\|^2 +\|M(s)\|^2%
)dv\bigg\}^{\frac{1}{2}}.
\end{equation}
\end{theorem}
\textit{Proof. } Under the hypothesis of the theorem,  Lemma
 \ref{lemaceg} implies
\begin{equation}\label{2.3.54}
E(t)=\int _0^t \dd\int_\Omega(\langle f,{u}_{,t}\rangle +\langle M,{P}_{,t}\rangle%
)dv ds,\  \forall \ t\geq 0.
\end{equation}%
By means of the Cauchy-Schwarz inequality we obtain%
\begin{equation}\label{2.3.55}
E(t)\leq \int _0^t\biggl\{\dd\int_\Omega(\| {u}_{,t}\|^2 +\|{P}_{,t}\|^2%
)dv\biggr\}^{\frac{1}{2}%
}g(s)ds,
\end{equation}%
for all $t\in[0,T].$  We define the function $\mathcal{Y}:[0,T]\rightarrow \mathbb{R}_+$\,
$
\mathcal{Y}(t)=[E(t)]^{\frac{1}{2}}.
$ This function is well defined because $E(\cdot)$ is a positive function. The inequality
(\ref{2.3.55}) becomes
\begin{equation}\label{2.3.58}
\mathcal{Y}^{2}(t)\leq \int _0^t\mathcal{Y}(s)g(s)ds,\
\forall\ t\geq 0.
\end{equation}%
By the Brezis' lemma given in Appendix (see \cite{VrabieDiff}, p. 47, Lemma 1.5.3)
 we deduce the
inequality
\begin{equation*}
\mathcal{Y}(t)\leq {\frac{1}{2}}\int _0^tg(s)ds,\ \forall \
t\geq 0,
\end{equation*}%
and the proof is complete.\hfill $\Box$

\subsection{Existence of the solution}\label{existences}

In this subsection, in order to establish an existence
theorem for the solution of the problem $(\mathcal{P})$  we use the results of
the semigroup theory of linear operators. First, we will rewrite the initial boundary value problem
$({\mathcal{P}})$ as an abstract Cauchy problem in a Hilbert space \cite{Pazy,Vrabie}. Let us define the space
\begin{equation}
\mathcal{X}\,{=}\,\big\{\,w=(u,v,P,K)\,|\,\ u{\in}{H}^1_0(\Omega),\quad v\in L^2(\Omega),\quad P{\in} H_0(\Curl; \Omega),\quad K\in L^2(\Omega)\big\}.
\end{equation}

On $\mathcal{X}$ we define the following bilinear form
\begin{align}\label{proscalar}
 (w_1,w_2)=\dd\int _\Omega\biggl(\langle v_1,v_2\rangle&+ \langle K_1,K_2\rangle+ \langle \C. \sym(\nabla u_1-P_1), \sym(\nabla u_2-P_2)\rangle
\vspace{2mm}\\&+ \langle \H.\sym P_1, \sym P_2\rangle+ \langle \L. \Curl P_1, \Curl P_2\rangle\biggr)dv,\notag
\end{align}
where
$w_1=(u_1,v_1,P_1,K_1)$ and $w_2=(u_2,v_2,P_2,K_2)$. Using the Lemma \ref{lemaii} and the same method as in the proof of Theorem \ref{cduid}, we observe that there is a positive constant $a_m$ such that
\begin{align}
a_m\big(\|v\|^2_{L^2(\Omega)}+ \|K\|^2_{L^2(\Omega)}+ \| \nabla u \|^2_{L^2(\Omega)}+ \|P\|^2_{L^2(\Omega)}+\|\Curl P \|^2_{L^2(\Omega)}\big)\leq (w,w),
\end{align}
where
$w=(u,v,K,P)\in \mathcal{X}$.

Hence, according with the symmetries \eqref{simetries}  we can conclude that the above bilinear form is an inner product on $\mathcal{X}$.

\begin{remark}
As in the proof of Theorem \ref{cduid} we have used that $h_m\neq 0$ in order to prove the above inequality. Hence, the above bilinear form $(\cdot,\cdot)$ is an inner product on $\mathcal{X}$ if $h_m\neq 0$.
\end{remark}

 Obviously, in view of \eqref{posdef}
 and of the following inequalities
 \begin{align}\label{inequs}
 &\|\sym(\nabla u-P)\|^2\leq 2\left( \|\sym \nabla u\|^2+\|\sym P\|^2\right),\notag\\
 &\|\sym \nabla u\|^2\leq \|\nabla u\|^2,\quad\quad
 \|\sym P\|^2\leq \|P\|^2,
 \end{align}
 we observe that there is also a positive constant $a_M$ such that
 \begin{align}
 (w,w)\leq a_M\big(\|v\|^2_{L^2(\Omega)}+ \|K\|^2_{L^2(\Omega)}+ \| \nabla u \|^2_{L^2(\Omega)}+ \|P\|^2_{L^2(\Omega)}+\|\Curl P \|^2_{L^2(\Omega)}\big)\, .
\end{align}

A direct consequence of the above  inequalities  is the fact that the norm induced by  $(\cdot,\cdot)$
is equivalent with the usual norm on  $\mathcal{X}$. Further, we introduce the operators
\begin{equation}
\hspace{-0.2cm}\dd\barr{crl} A_1\, w&=&v,\vspace{1.2mm}\\
A_2\, {w}&=&\dvg[ \C. \sym(\nabla u-P)],\\
A_3\, {w}&=&K,\vspace{1.2mm}\\
A_4\, {w}&=&- \crl[ \L.\crl\,P]+\C. \sym (\nabla u-P)-\H. \sym P\, ,\earr
\end{equation}
where all the derivatives  of the functions are understood in the sense of distributions. Let $\mathcal{A}$ be the operator
\begin{equation}\label{defopex}
\mathcal{A}=(A_1,A_2,A_3,A_4)
\end{equation}
with domain
\begin{equation}
\barr{crl}
\mathcal{D}(\mathcal{A})&=&\{{w}=(u,v,P,K)\in\mathcal{X}\ | \ \mathcal{A}{w}\in
\mathcal{X}\}.\earr
\end{equation}

We note that
${C}_0^\infty(\Omega)\!\times\!{C}_0^\infty(\Omega)\!\times\!{C}_0^\infty(\Omega)\!\times\!{C}_0^\infty(\Omega)$
is a dense subset of  $\mathcal{X}$  which is contained in $\mathcal{D}(\mathcal{A})$. Hence, $\mathcal{D}(\mathcal{A})$ is a dense subset of $\mathcal{X}$.

With the above definitions, the problem $({\mathcal{P}})$
can be transformed into the following abstract problem in the Hilbert
space
$\mathcal{X}$
\begin{equation}\label{prcauchy}
\frac{d{w}}{dt}(t)=\mathcal{A}{w}(t)+{\mathcal{F}}(t),\
\ {w}(0)={w}_0,
\end{equation}
where
\begin{equation}
{\mathcal{F}}(t)=\left({0},f,{0},M\right)
\end{equation} and
\begin{equation}
{w}_0=(u_0,\dot{u}_0, P_0, \dot{P}_0).
\end{equation}

%LEMMA 3.1

\begin{lemma} The operator $\mathcal{A}$ is dissipative\,\footnote{In fact $(\mathcal{A}{w},{w})= 0,
\text{ for all }{w} \in \mathcal{D}(\mathcal{A})$.}, i.e.
\begin{equation}\label{coerciv}
(\mathcal{A}{w},{w})\leq 0,
\text{ for all }{w} \in \mathcal{D}(\mathcal{A}) \ \ \text{in the inner product}\ \ (\cdot,\cdot) \ \ \text{defined in} \ \ \text{\eqref{proscalar}}.
\end{equation}
\end{lemma}
\textit{Proof}. Using the relations \eqref{formule} we find
that
\begin{align}
(\mathcal{A}{{w}},{w})=\dd\int _\Omega&\biggl(\langle {\rm Div} (\C. \sym(\nabla u-P)),v\rangle-\langle \crl(\L. \crl\, P),K\rangle\\&\notag+  \langle \C. \sym (\nabla u-P), K\rangle-\langle \H. \sym P,  K\rangle\\\notag &+ \langle \C. \sym(\nabla v-K), \sym(\nabla u-P)\rangle
\vspace{2mm}\notag\\&+\langle \H.\sym K,\sym P\rangle+ \langle \L. \Curl K, \Curl P\rangle\biggr)dv\notag\\\notag
=\dd\int _\Omega&\biggl(\sum\limits_{i=1}^3{\rm div} (v_{i} [\C. \sym(\nabla u-P)]_i)-\langle \C. \sym(\nabla u-P), \sym \nabla v\rangle\, \\\notag
&-\sum\limits_{i=1}^3{\rm div}\,  [(\L. \curl P)_i\times K_{i}]-\langle \L.\Curl\, P, \Curl\, K\rangle
\\& \notag+  \langle \C. \sym (\nabla u-P), K\rangle-\langle \H. \sym P,  K\rangle
\vspace{2mm}\notag\\
&\notag+ \langle \C. \sym(\nabla v-K), \sym(\nabla u-P)\rangle
\vspace{2mm}\notag\\&+ \langle \H.\sym K, \sym P\rangle+ \langle \L. \Curl K, \Curl P\rangle\biggr)dv\, .\notag
\end{align}

Hence, using the divergence theorem and the boundary conditions  $u=0$ and $P\times n=0$, we deduce
\begin{align}
(\mathcal{A}{{w}},{w})=\dd\int _\Omega\biggl(&-\langle \C. \sym(\nabla u-P), \sym(\nabla v-K)\rangle\,-\langle \L. \Curl\, P, \Curl\, K\rangle
\\& \notag-\langle \H. \sym  P,  \sym K\rangle+ \langle \C. \sym(\nabla v-K), \sym(\nabla u-P)\rangle
\vspace{2mm}\notag\\&+ \langle \H.\sym K, \sym P\rangle+ \langle \L. \Curl K, \Curl P\rangle\biggr)dv.\notag
\end{align}

The symmetries \eqref{simetries} assure that
\begin{align}
(\mathcal{A}{{w}},{w})=0, \ \ \text{for all} \ \ w\in \mathcal{D}(\mathcal{A})\, ,
\end{align}
 and the proof is complete.\hfill $\Box$

\begin{remark}
The dissipative condition $(\mathcal{A}{{w}},{w})\leq0, \ \ \text{for all} \ \ w\in \mathcal{D}(\mathcal{A})\,$ is already true for $h_m=0$ but this alone does not imply that $\mathcal{A}$ is dissipative, since for $h_m=0$ the bilinear form $(\cdot,\cdot)$ is not an inner product.
\end{remark}

\begin{lemma}\label{lemarange} The operator $\mathcal{A}$ satisfies the range condition, i.e.
\begin{equation}
{ \rm R}(I-\mathcal{A})=\mathcal{X}.
\end{equation}
\end{lemma}
\textit{Proof. }
 Let us consider
${w}^*=(u^*,v^*,P^*, K^*)\in \mathcal{X}$. We must show that
 the system
\begin{align}\label{exsis1}
\notag u-A_1{w}=&u^{*},\quad\quad\quad \ v-A_2{w}=v^{*}\,,\\
P-A_3{w}=&P^{*},\quad\quad \quad
K-A_4{w}=K^{*}
\end{align}
has a solution in  $\mathcal{D}(\mathcal{A})$.

By eliminating the functions   $v$ and $K$, we obtain for
the determination of the functions $u$ and $P$ the following system of equations
\begin{align}\label{exsis2}
L_1\mathbf{y}&\equiv u-\dvg[ \C. \sym(\nabla u-P)]=g_1,\\\notag
L_2\mathbf{y}&\equiv P+ \crl( \L.\crl\, P)-\C. \sym (\nabla u-P)+\H. \sym P=g_2\, .
\end{align}
where ${y}=(u,v,P,K)$,
\begin{equation}\label{g1g2}
 g_1=v^*+u^*, \quad\quad\quad g_2=K^*+P^*\, ,
\end{equation}
and all the derivatives  of the functions are understood in the sense of distributions.

We study this system in the following Hilbert space
\begin{equation}
\mathcal{Z}\!=\!{H}_0^1(\Omega)
\times{H}_0(\Curl;\Omega).
\end{equation}
We introduce the bilinear form
$\mathcal{B}:\mathcal{Z}\times\mathcal{Z}\rightarrow\mathbb{R}$
\begin{align}
\mathcal{B}({y},{\widetilde{y}})=\left\langle({L}_1{y},{L}_2{y}),
(\widetilde{u},\widetilde{P})\right\rangle_{{L}^2(\Omega)\times{L}^2(\Omega)}\, .
\end{align}
In view of relations \eqref{formule} and of the boundary conditions $u=0$ and $P\times n=0$, we have
\begin{align}
\mathcal{B}({y},{\widetilde{y}})=\dd\int _\Omega&\biggl(\langle u,\widetilde{u}\rangle-\langle \dvg (\C. \sym(\nabla u-P)),\widetilde{u}\rangle\\&\notag+\langle P,\widetilde{P}\rangle+\langle \H. \sym P,  \widetilde{P}\rangle+\langle \crl(\L. \crl\, P),\widetilde{P}\rangle  -\langle \C. \sym (\nabla u-P), \widetilde{P}\rangle\biggl)dv\\\notag
=\dd\int _\Omega&\biggl(\langle u,\widetilde{u}\rangle+\langle P,\widetilde{P}\rangle-\sum\limits_{i=1}^3{\rm div} (\widetilde{u}_{i} [\C. \sym(\nabla u-P)]_i)+\langle \C. \sym(\nabla u-P), \sym \nabla \widetilde{u}\rangle\, \\\notag
&+\langle \H. \sym P,  \widetilde{P}\rangle-\sum\limits_{i=1}^3{\rm div}\,  [(\L. \curl P)_i\times \widetilde{P}_{i}]+\langle \L. \Curl\, P, \Curl\, \widetilde{P}\rangle-  \langle \C. \sym (\nabla u-P), \widetilde{P}\rangle\biggl)dv
\\\notag
=\dd\int _\Omega&\biggl(\langle u,\widetilde{u}\rangle+\langle P,\widetilde{P}\rangle+\langle \C. \sym(\nabla u-P), \sym(\nabla \widetilde{u}-\widetilde{P})\rangle\,\\\notag&\ \ \ \ +\langle \H. \sym P, \sym \widetilde{P}\rangle +\langle \L. \Curl\, P, \Curl\, \widetilde{P}\rangle\biggl)dv\, ,
\end{align}
where
$\widetilde{y}=(\widetilde{u},\widetilde{P})$.
Let us define the linear operator  $l:\mathcal{Z}\rightarrow\mathbb{R}$
\begin{equation}\label{ling1g2}
l(\widetilde{y})=\left\langle({g}_1,{g}_2),(\widetilde{u},\widetilde{P})
\right\rangle_{{L}^2(\Omega)\times{L}^2(\Omega)}\, .
\end{equation}

From  the  Cauchy-Schwarz inequality and the Poincar\'{e} inequality we see that
the linear operator $l$ is bounded, i.e. there exists a positive constant $C$ such that
\begin{equation}
l({\widetilde{y}})\leq
C\|{\widetilde{y}}\|_\mathcal{Z}.
\end{equation}

Moreover, the Cauchy-Schwarz inequality leads us to
\begin{align}
{B}({y},\widetilde{{y}})\leq\dd\Bigg[&\int _\Omega\biggl(\| u\|^2+\| P\|^2+\langle \C. \sym(\nabla u-P), \sym(\nabla {u}-{P})\rangle\, \\\notag
&\ \ \ \ \ \ \ +\langle \H. \sym P, \sym {P}\rangle+\langle \L. \Curl\, P, \Curl\, {P}\rangle\biggl) dv\Bigg]^{\frac{1}{2}}
\\\notag
&\times \Bigg[\int _\Omega\biggl(\| \widetilde{u}\|^2+\| \widetilde{P}\|^2+\langle \C. \sym(\nabla \widetilde{u}-\widetilde{P}), \sym(\nabla {\widetilde{u}}-{\widetilde{P}})\rangle\,\\\notag&\ \ \ \ \ \ \  +\langle \H. \sym \widetilde{P}, \sym {\widetilde{P}}\rangle+\langle \L. \Curl\, \widetilde{P}, \Curl\, {\widetilde{P}}\rangle\biggl)dv \Bigg]^{\frac{1}{2}}\, .
\end{align}

In view of \eqref{posdef} we obtain that
\begin{align}
\mathcal{B}({y},\widetilde{y})\leq \dd \,C\, \Bigg[&\int _\Omega\biggl(\| u\|^2+\| P\|^2+\| \sym(\nabla u-P)\|^2+\| \sym P\|+\|\Curl\, P\|^2\biggl)dv\Bigg]^{\frac{1}{2}}
\\\notag
&\times \Bigg[\int _\Omega\biggl(\| \widetilde{u}\|^2+\| \widetilde{P}\|^2+\| \sym(\nabla \widetilde{u}-\widetilde{P})\|^2+\| \sym \widetilde{P}\|+\|\Curl\, \widetilde{P}\|^2\biggl)dv\Bigg]^{\frac{1}{2}}\, ,
\end{align}
where $C$ is a positive constant.

Hence, using the inequalities \eqref{inequs} and the Poincar\'{e} inequality, we can find a positive constant $C$ such that
\begin{align}
\mathcal{B}({y},\widetilde{y})\leq \dd \,C\, \Bigg[&\int _\Omega\biggl(\| u\|^2+\| \nabla u\|^2+\| P\|^2+\|\Curl\, P\|^2\biggl)dv\Bigg]^{\frac{1}{2}}
\\\notag
&\times \Bigg[\int _\Omega\biggl(\| \widetilde{u}\|^2+\| \nabla \widetilde{u}\|^2+\| \widetilde{P}\|^2+\|\Curl\, \widetilde{P}\|^2\biggl)dv\Bigg]^{\frac{1}{2}}\leq \dd \,C\, \|y\|_{\mathcal{Z}}\|\widetilde{y}\|_{\mathcal{Z}}\, ,
\end{align}
which means that $\mathcal{B}$ is bounded. On the other hand, we have
\begin{align}
\mathcal{B}({y},{y})
=\dd\int _\Omega\biggl(&\| u\|^2+\| P\|^2+\langle \C. \sym(\nabla u-P), \sym(\nabla {u}-{P})\rangle\,\\& +\langle \H. \sym P, \sym {P}\rangle+\langle \L. \Curl\, P, \Curl\, {P}\rangle\biggl)\, dv\,\notag
\end{align}
for all
${y}=({u},{P})\in \mathcal{Z}$. The inequality \eqref{ineqmed} shows that
\begin{align}
\langle \C. \sym(\nabla u-P), \sym(\nabla u-P)\rangle&\geq {c_m}\bigg[(1-\delta)\| \sym \nabla u\|^2+\bigg(1 -\frac{1}{\delta}\bigg)\|\sym P\|^2\bigg]\, ,
\end{align}
for all $\delta>0$. Hence, we deduce that
\begin{align}
 {c_m}(1-\delta)\| \sym \nabla u \|^2+\bigg({c_m} +1-\frac{{c_m}}{\delta}\bigg)\|\sym P\|^2 &\leq \langle \C. \sym(\nabla u-P), \sym(\nabla u-P)\rangle+ \|\,\sym P\|^2\notag\\
& \leq \langle \C. \sym(\nabla u-P), \sym(\nabla u-P)\rangle+ \|\,P\|^2
\notag\\
& \leq {c_M}\| \sym(\nabla u-P)\|^2+ \|\,P\|^2.
\end{align}

If we choose $\delta$ so that
$$
\frac{{c_m} }{{c_m} +1}<\delta<1
$$
we obtain a positive constant $C_1$ so that
\begin{align}\label{estbun}
C_1\big( \| \sym \nabla u\|^2+\|\sym P\|^2\big)\leq \| \sym(\nabla u-P)\|^2+ \|\,P\|^2\, .
\end{align}

Moreover, as a consequence of  the assumptions \eqref{posdef} we have
\begin{align}
\mathcal{B}({y},{y})
\geq \dd\int _\Omega\biggl(&\| u\|^2+\| P\|^2+{c_m}\| \sym(\nabla u-P)\|^2+h_m\| \sym P\|^2+{(L_c)}_m \|\Curl P\|^2\biggl)\, dv\,\\\notag
\geq \dd\int _\Omega\biggl(&\| u\|^2+\min\{{c_m},1\}(\| P\|^2+\| \sym(\nabla u-P)\|^2)+h_m\| \sym P\|^2+{(L_c)}_m \|\Curl P\|^2\biggl)\, dv\,
\\\notag
\geq \dd\int _\Omega\biggl(&\| u\|^2+C_1\min\{{c_m},1\}\big( \| \sym \nabla u\|^2+\|\sym P\|^2\big)+h_m\| \sym P\|^2+{(L_c)}_m \|\Curl P\|^2\biggl)\, dv\,
\\\notag
\geq \dd\int _\Omega\biggl(&\min\{1,C_1\min\{{c_m},1\},{(L_c)}_m\}\big(\| u\|^2+ \| \sym\nabla u\|^2+\|\sym P\|^2+\|\Curl P\|^2\big)+h_m\| \sym P\|^2\biggl)\, dv\, .
\end{align}

Using \eqref{Korn} and \eqref{Neff2} we deduce
\begin{align}\label{oi}
\mathcal{B}({y},{y})
\geq \dd\int _\Omega\biggl(&C\big(\| u\|^2+ \| \nabla u\|^2+\|P\|^2+\|\Curl P\|^2\big)+h_m\| \sym P\|^2\biggl)\, dv\,
\geq \dd C\|y\|^2_{\mathcal{Z}}\, ,
\end{align}
where $C$ is a positive constant. Hence $\mathcal{B}(\cdot,\cdot)$ is coercive.

Using the Lax-Milgram theorem we prove the existence of a solution of the system (\ref{exsis2}) in $\mathcal{Z}$, i.e.
 \begin{equation}
 u\in {H}_0^1(\Omega), \quad\quad\quad \ P\in{H}_0(\Curl;\Omega)\ .
 \end{equation}

 Hence, $v$ and $K$ will be given by
  \begin{equation}
 v=u-u^*,\quad\quad\quad  K=P-P^*\ .
 \end{equation}

 Moreover, because $w^*\in \mathcal{X}$ it follows that $u^*\in H_0^1(\Omega)$ and $P^*\in H_0(\Curl;\Omega)$. Thus
 \begin{equation}
 v\in {H}_0^1(\Omega)\subset L^2(\Omega), \quad\quad\quad  \ K\in{H}_0(\Curl;\Omega)\subset L^2(\Omega)\,.
 \end{equation}

 Let us remark that
  \begin{equation}
\mathcal{D}(\mathcal{A})\subset  {H}_0^1(\Omega)\times  {H}_0^1(\Omega)\times{H}_0(\Curl;\Omega)\times {H}_0(\Curl;\Omega)\subset \mathcal{X} ,
 \end{equation}
 and that, until now, we have proved that the equation
 \begin{align}
 w-\mathcal{A}\,w=w^*
 \end{align}
 has a solution
 \begin{equation}
 w=(u,v,P,K)\in  {H}_0^1(\Omega)\times  {H}_0^1(\Omega)\times{H}_0(\Curl;\Omega)\times {H}_0(\Curl;\Omega)\subset \mathcal{X}
 \end{equation}
 for all $w^*\in \mathcal{X}$. Hence, $-\mathcal{A}\,w=w^* -w\in  \mathcal{X}$, which  shows that $w\in \mathcal{D}(\mathcal{A})$, and
$ w-\mathcal{A}\,w=w^* $. This implies that we have the desired solution
of system (\ref{exsis1}) and the proof is complete.\hfill $\Box$

\begin{remark}
In the proof of the range condition from the previous lemma, we have not used that $h_m>0$. The result holds true if $\H=0$. However, the bilinear form $(\cdot,\cdot)$ is an inner product if $h_m>0$ and we can not prove this fact if $\H=0$. Our existence result needs that $\mathcal{A}$ is dissipative in the inner product $(\cdot,\cdot)$. Hence, the existence result below holds true only for $h_m>0$.
\end{remark}

%THEOREM 3.1

\begin{theorem}
The operator  $\mathcal{A}$ defined by  {\rm (\ref{defopex})} generates a  $C_0$-contractive semigroup
in $\mathcal{X}$.
\end{theorem}

{\it Proof.}
 The proof follows using the previous lemmas and the
Lumer--Phillips corollary
to the Hille--Yosida theorem \cite{Pazy} given in  Appendix. \hfill $\Box$

%THEOREM 3.2

\begin{theorem}\label{teorexs} Assume that   $f,M\, \in C^1([0,T);L^2(\Omega))$,
${w}_0\in\mathcal{D}(\mathcal{A})$ and the fourth order elasticity tensors $\C$, $\L$ and $\H$ are positive definite and satisfy the symmetries \eqref{simetries}. Then, there exists a unique solution  $${w}\!\in\!
{C}^1((0,T);\mathcal{X})\cap
{C}^0([0,T);\mathcal{D}(\mathcal{A}))$$ of the Cauchy problem
{\rm (\ref{prcauchy}).} \hspace{12cm}
\end{theorem}

{\it Proof.}
The proof follows from the results concerning the abstract
Cauchy problem from the Appendix (see \cite{Pazy,Vrabie}). \hfill $\Box$

%COROLLARY 3.2

\begin{corollary} In the hypothesis of Theorem  {\rm \ref{teorexs}} we have the following estimate
\begin{equation}\barr{crl}
&&\|{w}(t)\|_\mathcal{X}\leq
\|{w}_0(t)\|_\mathcal{X}+C \dd\int
_0^{T}\left(\|f(s)\|_{{L}^2(\Omega)}+\|M(s)\|_{{L}^2(\Omega)}\right)ds , \text{\quad for all\quad } t\in[0,T],\earr
\end{equation}
where $C$ is a positive constant.
\hspace{9.5cm}
\end{corollary}

{\it Proof.}  For the proof of this Corollary we use the fact that
the semigroup  generated by $\mathcal{A}$ is contractive and apply the Duhamel Principle (see the Appendix). \hfill $\Box$

\section{Another further relaxed  problem}
\setcounter{equation}{0}

In this section, we weaken our energy expression further in the following model, where the corresponding elastic energy depends now only on the set of  {\it independent constitutive variables}
 \begin{align}
 \varepsilon_e=\sym(\nabla u-P),\quad\quad\quad \dev \varepsilon_p=\dev\sym P, \quad\quad\quad \dev\alpha=-\dev \Curl P.
 \end{align}
  In this model, it is neither implied that $P$ remains symmetric, nor that $P$ is trace-free, but only the trace free symmetric part of the micro-distortion $P$ and the trace-free part of the micro-dislocation tensor $\alpha$ contribute to the stored energy.

\subsection{Formulation of the problem}
The model in its general anisotropic form is:
\begin{align}\label{eqdev}
u_{,tt}&={\dvg}[ \C. \sym(\nabla u-P)]+f\, ,\\\notag
P_{,tt}&=- {\Curl}[ \dev [\L.\dev \Curl P]]+\C. \sym (\nabla u-P)-\H. \dev \sym P+M\, \ \ \ \text {in}\ \ \  \Omega\times [0,T].
\end{align}

In the isotropic case the model becomes
\begin{align}\label{eqisoup2}
 u_{,tt}&=\dvg[2\mu_e \sym(\nabla u-P)+\lambda_e \tr(\nabla u-P){\cdp} \id]+f\, ,\\\notag
P_{,tt}&=-\Curl [\alpha_1 \dev\sym \Curl P+\alpha_2 \skew \Curl P ]\\&\quad\ +2\mu_e \sym(\nabla u-P)+\lambda_e \tr(\nabla u-P){\cdp} \id-2\mu_h \dev\sym P+M\,  \ \ \ \text {in}\ \ \  \Omega\times [0,T].\notag
\end{align}

To the system of partial differential equations of this model we adjoin the weaker boundary conditions
\begin{align} \label{bcdev}
{u}({x},t)=0, \ \ \ \quad \quad {P}_i({x},t)\times n(x) =0, \ \ \ i=1,2,3, \ \ \
\ ({x},t)\in\partial \Omega\times [0,T],
\end{align}%
and the nonzero initial conditions
\begin{align}\label{icdev}
&{u}({x},0)={u}_0(
x), \quad\quad\quad \dot{u}({x},0)=\dot{u}_0(
x),\quad\quad\quad {P}({x},0)={P}_0(
x), \quad\quad\quad  \dot{P}({x},0)=\dot{P}_0(
x),\ \ \text{\ \ }{x}\in \bar{\Omega},
\end{align}%
where  ${u}_0, \dot{u}_0, {P}_0$ and $\dot{P}_0$ are prescribed functions.

We remark again that $P$ is not trace-free in this formulation and no projection is performed. We denote the new problem defined by the above equations, the boundary conditions \eqref{bcdev} and the initial conditions \eqref{icdev} by $(\widetilde{\mathcal{P}})$.

\subsection{Energy conservation}
 To a solution   $\{u,P\}$ of the problem  $(\widetilde{\mathcal{%
P}})$ we associate the
total energy
\begin{align}
2\,\widetilde{E}(t)=&\int_\Omega\bigg(\|u_{,t}\|^2+ \|P_{,t}\|^2+ \langle \C.\, \sym(\nabla u-P), \sym(\nabla u-P)\rangle
\vspace{2mm}\notag\\& \ \ \ \ \ \ \ \ \ \ + \langle \H.\,\dev \sym P, \dev \sym P\rangle + \langle \L.\, \dev \Curl P, \dev \Curl P\rangle\bigg)dv\, .
\end{align}
We observe that since $\mathbb{H}$ is positive definite on $\Sym(3)$,  in view of \eqref{posdef} we also have the estimate
\begin{align}
h_m\|\dev \sym P \|^2\leq &\langle \H. \dev \sym P,\dev \sym P\rangle\leq h_M\| \dev \sym P\|^2\,  \ \ \text{for all }\ \ P\in \mathbb{R}^{3\times3}\,.
\end{align}

It is easy to see that taking in \eqref{formule}
\begin{align}
R_i=[\dev (\L. \dev\Curl P)]_i, \quad \  S_i=P_i\, ,
\end{align}
we have
\begin{align}
\sum\limits_{i=1}^3{\rm div}\,  ([\dev (\L. \dev\Curl P)]_i\times P_{i,t})=&\sum\limits_{i=1}^3\langle P_{i,t}, \curl \, [\dev (\L. \dev\Curl P)]_i\rangle\\&-\sum\limits_{i=1}^3\langle [\dev (\L. \dev\Curl P)]_i, \curl\,  P_{i,t}\rangle\, .\notag
\end{align}

Thus, in terms of tensors, we have the following equality
\begin{align}\label{fcurlm}
\langle P_{,t},\Curl[ \, \dev(\L. \dev\Curl \,  P)]\rangle=\langle \dev &[\L. \dev (\Curl\, P)], \Curl\, P_{,t}\rangle\\
\notag&+\sum\limits_{i=1}^3{\rm div}\,  ([\dev(\L. \dev\Curl P)]_i\times P_{i,t})
\, .
\end{align}

Moreover, using that
\begin{equation}\label{ABdev}
\langle \dev A,B\rangle=\langle A, \dev B\rangle, \ \  \text{for all} \ \ \ A,B\in \mathbb{R}^{3\times 3},
\end{equation}
we have
\begin{align}\label{fcurlm}
\langle P_{,t},\Curl[ \, \dev(\L. \dev(\Curl \,  P))]\rangle=\langle  \L. \dev& \Curl\, P, \dev \Curl\, P_{,t}\rangle
\\\notag
&+\sum\limits_{i=1}^3{\rm div}\,  ([\dev(\L. \dev\Curl P)]_i\times P_{i,t})\, .
\end{align}

In view of \eqref{eqdev} and  using the symmetries \eqref{simetries}, we get
\begin{align}
\frac{1}{2}\frac{\partial }{\partial t}\bigg(&\|u_{,t}\|^2+ \|P_{,t}\|^2+ \langle \C. \sym(\nabla u-P), \sym(\nabla u-P)\rangle
\vspace{2mm}\notag\\&+ \langle \H.\dev \sym P, \dev \sym P\rangle+ \langle \L. \dev \Curl P, \dev \Curl P\rangle\, \bigg) \\\notag
&=\sum\limits_{i=1}^3 \langle {\rm div}([\C. \sym(\nabla u-P)]_iu_{i,t})+\sum\limits_{i=1}^3{\rm div}( P_{i,t}\times [\dev(\L. \dev \Curl P)]_i)+\langle f,u_{,t}\rangle+\langle M,P_{,t}\rangle\, .\notag
\end{align}

Using the divergence theorem and the boundary conditions $u=0$ and $P\times n=0$, from the above identity we can conclude:
\begin{remark}
{\rm (Energy conservation for the problem $(\widetilde{\mathcal{P}})$)}
 If $\{u,P\}$
 is a solution of the problem $(\widetilde{\mathcal{P}})$ corresponding to the loads  $\mathcal{I}=\{f,M, {u}_0, \dot{u}_0, {P}_0, \dot{P}_0\}$, then  we have
\begin{equation}\label{consdev}
\widetilde{E}(t)=\widetilde{E}(0)+\int_0^t\, \Pi(s)ds\, ,
\end{equation}
for every time $t\in [0,T]$.
\end{remark}

\subsection{Uniqueness and continuous dependence of the solution}

In order to prove the uniqueness and the continuous dependence of the solution with respect to given data, we will use the estimate given by the following lemma:
\begin{lemma}\label{esdevlemma}
For all $u\in H^1(\Omega)$ and $P\in H(\Curl;\Omega)$, the following estimate holds true
\begin{align}\label{estdev}
a\bigg(\| \nabla u\|^2_{L^2(\Omega)}&+\|\dev \sym P\|^2_{L^2(\Omega)}\bigg)\\\notag&\leq  \int_{\Omega}\bigg(\langle\C. \sym(\nabla u-P), \sym(\nabla u-P)\rangle+ \langle \H.\dev \sym P, \dev \sym P\rangle\bigg)dv\, ,
\end{align}
 where   $a$ is a positive constant.
\end{lemma}

{\it Proof.} First, we note that
\begin{align}\label{devA}
\|X\|^2=\| \dev X\|^2+\dd\frac{1}{3}(\tr X)^2,
\end{align}
for all $X\in\mathbb{R}^{3\times 3}$. Using this identity and the properties \eqref{posdef} we deduce
\begin{align}\label{ineqmed0}
\langle \C. \sym(\nabla u-P), \sym(\nabla u-P)\rangle &+ \langle \H.\dev \sym P, \dev \sym P\rangle \\\notag&\geq {c_m} \|\sym(\nabla u-P)\|^2+ h_m\|\dev \sym P\|^2\notag\\
&\geq {c_m} \|\dev\sym(\nabla u-P)\|^2+ h_m\|\dev \sym P\|^2\notag\\
&\geq {c_m}(1-\delta)\|\dev \sym \nabla u\|^2+\bigg({c_m} -\frac{{c_m}}{\delta}+h_m\bigg)\|\dev\sym P\|^2\notag
\end{align}
for all $\delta>0$. If we choose $\delta$ so that
$$
\frac{{c_m} }{{c_m} +h_m}<\delta<1
$$
we have that there is a positive constant $a_1$ so that
\begin{align}\hspace*{-7mm}
\langle \C. \sym(\nabla u-P), \sym(\nabla u-P)\rangle+ \langle \H.\dev\sym P,& \dev\sym P\rangle\\\notag
&\geq a_1\big( \| \dev\sym \nabla u\|^2+\|\dev\sym P\|^2\big)\,.
\end{align}
The estimate \eqref{estdev} follows from Theorem \ref{BNPS4}. \hfill $\Box$

\begin{theorem}\label{cduiddev}{\rm (Continuous dependence upon initial data)}  Let $\{u,P\}$ %
be a solution of the problem $(\widetilde{\mathcal{P}})\ $
with the external data system $\mathcal{I}=\{0,0, {u}_0, \dot{u}_0, {P}_0, \dot{P}_0\}$.
 Then, there is a positive constant $a$ so that
\begin{equation}\label{cdidev}
a\big(\|u_{,t}\|^2_{L^2(\Omega)}+ \|P_{,t}\|^2_{L^2(\Omega)}+ \| \nabla u\|^2_{L^2(\Omega)}+\|P\|^2_{L^2(\Omega)}+\|\Curl P\|^2_{L^2(\Omega)}\big)\leq E(0), \quad \ \text{for all} \ \  t\in[0,T]\, .
\end{equation}
\end{theorem}

{\it Proof.} Using the  conservation law \eqref{consdev} and  the estimate \eqref{estdev} we have
 that there is a positive constant $a$ so that
\begin{align}
a\big(\|u_{,t}\|^2_{L^2(\Omega)}+ \|P_{,t}\|^2_{L^2(\Omega)}+ \| \nabla u\|^2_{L^2(\Omega)}+\|\dev \sym P\|^2_{L^2(\Omega)}+ \|\dev \Curl P \|^2_{L^2(\Omega)}\big)\leq E(0), \quad \ \text{for all} \ \  t\in[0,T]\, .\notag
\end{align}

Because  $P\in{\rm H}_0({\rm Curl}\, ; \Omega)$ and $u\in H_0^1(\Omega)$, in view of \eqref{BNPS}, \eqref{BNPSe3}, there are the positive constants  $C_{DSDC}$ and $C_{DD}$ \cite{NeffPaulyWitsch,NPW2,NPW3,BNPS1,BNPS2}, such that
\begin{align}\label{BNPSagain}
\| P\|_{L^2(\Omega)}+\| \Curl P\|_{L^2(\Omega)}\leq (C_{DSDC}+C_{DD})\left(\|\dev \sym P\|_{L^2(\Omega)}+\|\dev \Curl P\|_{L^2(\Omega)}\right)\, .
\end{align}

Hence, we can find a positive constant $a$ so that the inequality \eqref{cdidev} is satisfied.

\begin{corollary} {(\rm Uniqueness)}
Any two solutions of the problem  $(\widetilde{\mathcal{P}})$ are equal.%
\end{corollary}

For the modified problem $(\widetilde{\mathcal{P}})$ we also have  the continuous data dependence of solution upon the
supply terms  $\{f,M\}$.
\begin{remark}
{\rm (Continuous dependence upon the supply terms)} %
 Let  $\{u,P\}$ be a solution of the problem $(\widetilde{\mathcal{P}})$
corresponding to the external data system
$\mathcal{I}=\{f,M,0,0,0,0\}.$
Then, for all $t\in I$ we have
\begin{equation}
\sqrt{a}\big(\|u_{,t}\|^2_{L^2(\Omega)}+ \|P_{,t}\|^2_{L^2(\Omega)}+ \| \nabla u\|^2_{L^2(\Omega)}+\|P\|^2_{L^2(\Omega)}+\|\Curl P\|^2_{L^2(\Omega)}\big)^{\frac{1}{2}}\leq \frac{1}{2}\int _0^t g(s)ds.
\end{equation}%
\end{remark}
\textit{Proof. } The proof of this remark is similar to the proof of Theorem \ref{cdst}.\hfill $\Box$
\subsection{Existence of the solution}

In this subsection, we study the existence of  the solution of the problem $(\widetilde{\mathcal{P}})$. Because the method is similar with that used in Section \ref{existences} we only point out the differences which arise for our modified problem.

We consider the same Hilbert space $\mathcal{X}$ as defined in Section \ref{existences} and we define the following bilinear form
\begin{align}\label{proscalar}
 ((w_1,w_2))=\dd\int _\Omega\biggl(\langle v_1,v_2\rangle&+ \langle K_1,K_2\rangle+ \langle \C. \sym(\nabla u_1-P_1), \sym(\nabla u_2-P_2)\rangle
\vspace{2mm}\notag\\&+ \langle \H.\dev\sym P_1, \dev\sym P_2\rangle+ \langle \L. \dev \Curl P_1, \dev \Curl P_2\rangle\biggr)dv,\notag
\end{align}
where
$w_1=(u_1,v_1,P_1,K_1)$ and $w_2=(u_2,v_2,P_2,K_2)$. Let us remark that in view of \eqref{estdev}, there is a positive constant $a$ such that
\begin{align}
a\big(\|v\|^2_{L^2(\Omega)}+ \|K\|^2_{L^2(\Omega)}+ \| \nabla u \|^2_{L^2(\Omega)}+ \|\dev \sym P\|^2_{L^2(\Omega)}+\|\dev \Curl P \|^2_{L^2(\Omega)}\big)\leq ((w,w)),
\end{align}
where
$w=(u,v,K,P)\in \mathcal{X}$. In other words, using Corollary \ref{BNPS3} we have
\begin{align}
a\big(\|v\|^2_{L^2(\Omega)}+ \|K\|^2_{L^2(\Omega)}+ \| \nabla u \|^2_{L^2(\Omega)}+ \| P\|^2_{L^2(\Omega)}+\| \Curl P\|^2_{L^2(\Omega)}\big)\leq ((w,w))\, ,
\end{align}
which implies that the  bilinear form $((\cdot,\cdot))$ is an inner product on $\mathcal{X}$. We also have that there is  a positive constant $C$ such that
 \begin{align}
 ((w,w))\leq C\big(\|v\|^2_{L^2(\Omega)}+ \|K\|^2_{L^2(\Omega)}+ \| \nabla u \|^2_{L^2(\Omega)}+ \| P\|^2_{L^2(\Omega)}+ \|\dev \sym P\|^2_{L^2(\Omega)}+\|\dev \Curl P \|^2_{L^2(\Omega)}\big)\, .\notag
\end{align}

Hence, using Theorem \ref{BNPS} we obtain
 \begin{align}
 ((w,w))\leq C\big(\|v\|^2_{L^2(\Omega)}+ \|K\|^2_{L^2(\Omega)}+ \| \nabla u \|^2_{L^2(\Omega)}+ \| P\|^2_{L^2(\Omega)}+\|\Curl P\|^2_{L^2(\Omega)}\big)\, .
\end{align}

Thus,  the norm induced by  $((\cdot,\cdot))$
is equivalent with the usual  norm   on $\mathcal{X}$. We consider  the operators
\begin{equation}
\hspace{-0.2cm}\dd\barr{crl} \widetilde{A}_1\, w&=&v,\vspace{1.2mm}\\
\widetilde{A}_2\, {w}&=&\dvg[ \C. \sym(\nabla u-P)],\\
\widetilde{A}_3\, {w}&=&K,\vspace{1.2mm}\\
\widetilde{A}_4\, {w}&=&- \crl[ \dev [\L.\dev\crl\, P]]+\C. \sym (\nabla u-P)-\H. \dev \sym P\, ,\earr
\end{equation}
where all the derivatives  of the functions are understood in the sense of distributions, and the operator $\widetilde{\mathcal{A}}$
\begin{equation}\label{defopex}
\widetilde{\mathcal{A}}=(\widetilde{A}_1,\widetilde{A}_2,\widetilde{A}_3,\widetilde{A}_4)
\end{equation}
with the domain
\begin{equation}
\barr{crl}
\mathcal{D}(\widetilde{\mathcal{A}})&=&\{{w}=(u,v,P,K)\in\mathcal{X} \ | \ \widetilde{\mathcal{A}}{w}\in
\mathcal{X}\}.\earr
\end{equation}

A similar method with that considered in Section \ref{existences} gives us the following existence result:

\begin{theorem}\label{teorexs} Assume that   $f,M\, \in C^1([0,T);L^2(\Omega))$,
${w}_0\in\mathcal{D}(\widetilde{\mathcal{A}})$ and the fourth order elasticity tensors $\C$, $\L$ and $\H$ are positive definite and satisfy the symmetries \eqref{simetries}. Then, there exists a unique solution
$${w}\!\in\!
{C}^1((0,T);\mathcal{X})\cap
{C}^0([0,T);\mathcal{D}(\widetilde{\mathcal{A}}))
$$
 of the following Cauchy problem
\begin{equation}\label{prcauchy}
\frac{d{w}}{dt}(t)=\widetilde{\mathcal{A}}{w}(t)+{\mathcal{F}}(t),\
\ {w}(0)={w}_0,
\end{equation}
where
\begin{equation}
{\mathcal{F}}(t)=\left({0},f,{0},M\right)
\end{equation} and
\begin{equation}
{w}_0=(u_0,\dot{u}_0, P_0, \dot{P}_0).
\end{equation}
Moreover, we have the estimate
\begin{equation}\barr{crl}
&&\|{w}(t)\|_\mathcal{X}\leq
\|{w}_0(t)\|_\mathcal{X}+C \dd\int
_0^{T}\left(\|f(s)\|_{{L}^2(\Omega)}+\|M(s)\|_{{L}^2(\Omega)}\right)ds ,\ \text{\quad for all \quad} t\in[0,T],\earr
\end{equation}
where $C$ is a positive constant.
\end{theorem}

{\it Proof.} It is easy to prove that the operator $\widetilde{\mathcal{A}}$ is dissipative.  In the following we prove that the operator $\widetilde{\mathcal{A}}$ satisfies the range condition
\begin{equation}
{ \rm R}(I-\widetilde{\mathcal{A}})=\mathcal{X}.
\end{equation}

 Let us define the operator $\widetilde{L}_2:\mathcal{X}\rightarrow\mathcal{X}$ by
\begin{align}\label{exsis3}
\widetilde{L}_2\mathbf{y}&\equiv P+ \crl[ \dev[\L.\dev \crl\, P]]+\C. \sym (\nabla u-P)-\H.\dev  \sym P\, ,
\end{align}
where ${y}=(u,v,P,K)\in\mathcal{X}$ and all the derivatives  of the functions are understood in the sense of distributions. We consider the Hilbert space
\begin{equation}
\mathcal{Z}\!=\!{H}_0^1(\Omega)
\times{H}_0(\Curl;\Omega)\, .
\end{equation}

On $\mathcal{Z}$ we consider the bilinear form
$\widetilde{\mathcal{B}}:\mathcal{Z}\times\mathcal{Z}\rightarrow\mathbb{R}$
\begin{align}
\mathcal{B}({y},{\widetilde{y}})=\left\langle({L}_1{y},\widetilde{L}_2{y}),
(\widetilde{u},\widetilde{P})\right\rangle_{{L}^2(\Omega)\times{L}^2(\Omega)},
\end{align}
where $L_1$ is given by \eqref{exsis2}, and the linear bounded operator  $l:\mathcal{Z}\rightarrow\mathbb{R}$
is given by \eqref{ling1g2}. In view of  the boundary conditions we have
\begin{align}
\widetilde{\mathcal{B}}({y},{\widetilde{y}})
=\dd\int _\Omega&\biggl(\langle u,\widetilde{u}\rangle+\langle P,\widetilde{P}\rangle+\langle \C. \sym(\nabla u-P), \sym(\nabla \widetilde{u}-\widetilde{P})\rangle\,\\\notag&\ \ \ \ +\langle \H.\dev \sym P,\dev \sym \widetilde{P}\rangle+\langle \L.\dev \Curl\, P,\dev \Curl\, \widetilde{P}\rangle\biggl)dv\, ,
\end{align}
where
$\widetilde{y}=(\widetilde{u},\widetilde{P})$. The  Cauchy-Schwarz inequality, the Poincar\'{e} inequality and the relations \eqref{posdef}, \eqref{inequs}   and \eqref{devA} lead us to the estimate
\begin{align}
\widetilde{\mathcal{B}}({y},\widetilde{y})
\leq \dd \,C\, \Bigg[&\int _\Omega\biggl(\| u\|^2+\| \nabla u\|^2+\| P\|^2+\| \dev P\|^2+\|\dev\Curl\, P\|^2\biggl)dv\Bigg]^{\frac{1}{2}}
\\\notag
&\times \Bigg[\int _\Omega\biggl(\| \widetilde{u}\|^2+\| \nabla \widetilde{u}\|^2+\| \widetilde{P}\|^2+\|\dev \widetilde{P}\|^2+\|\dev\Curl\, \widetilde{P}\|^2\biggl)dv\Bigg]^{\frac{1}{2}}\\\notag\leq \dd \,C\, \Bigg[&\int _\Omega\biggl(\| u\|^2+\| \nabla u\|^2+\| P\|^2+\|\Curl\, P\|^2\biggl)dv\Bigg]^{\frac{1}{2}}
 \Bigg[\int _\Omega\biggl(\| \widetilde{u}\|^2+\| \nabla \widetilde{u}\|^2+\| \widetilde{P}\|^2+\|\Curl\, \widetilde{P}\|^2\biggl)dv\Bigg]^{\frac{1}{2}}\\\notag
&\!\!\!\!\!\!\!\!\!\!\!\!\!\!\!\!\!\!\leq \dd \,C\, \|y\|_{\mathcal{Z}}\|\widetilde{y}\|_{\mathcal{Z}}\, ,
\end{align}
where $C$ is a positive constant. This means that $\widetilde{\mathcal{B}}$ is bounded. On the other hand, we have
\begin{align}
\widetilde{\mathcal{B}}({y},{y})
=\dd\int _\Omega\biggl(&\| u\|^2+\| P\|^2+\langle \C. \sym(\nabla u-P), \sym(\nabla {u}-{P})\rangle\,\\& +\langle \H. \dev\sym P, \dev\sym {P}\rangle+\langle \L. \dev\Curl\, P,\dev \Curl\, {P}\rangle\biggl)\, dv\,\notag
\end{align}
for all
${y}=({u},{P})\in \mathcal{Z}$. Moreover, as a consequence of  the assumptions \eqref{posdef} and of estimate \eqref{estdev}, we can find a new positive constant $C$ so that
\begin{align}\label{oi0}
\mathcal{B}({y},{y})
\geq C\dd\int _\Omega\biggl(&\| u\|^2+ \| \nabla u\|^2+\|\dev \sym P\|^2+\|\dev \Curl P\|^2\biggr)\, dv\, .
\end{align}

Now, using the Corollary \ref{BNPS3} we deduce that $\mathcal{B}(\cdot,\cdot)$ is coercive. Using the Lax-Milgram theorem and similar arguments as in the proof of Lemma \ref{lemarange} it follows that the operator $\widetilde{\mathcal{A}}$ satisfies the range condition. Hence, the hypothesis of the general existence theorem for the abstract
Cauchy problem are satisfied and the proof is complete.\hfill $\Box$
\begin{figure}
\setlength{\unitlength}{1mm}
\begin{center}
\begin{picture}(10,220)
\thicklines
\put(-40,225){\footnotesize{\bf The family of relaxed micromorphic dislocation models }}
\put(0,198){\oval(51,41)}
\put(-20,215){\footnotesize{\bf Relaxed micromorphic }}
\put(-20,211){\footnotesize{\bf dislocation}}
\put(-20,207){\footnotesize{12 dof $(u,P)$, well-posed$^*$}}
\put(-20,203){\footnotesize{$\sigma$ symmetric,\quad $\sigma={\mathbb{C}}.\, \varepsilon_e$}}
\put(-20,199){\footnotesize{isotropic: 6+3 parameters}}
\put(-20,195){\footnotesize{no coupling: 4+3 parameters}}
\put(-20,191){\footnotesize{\it constitutive variables:}}
\put(-20,187){\footnotesize{$\varepsilon_e=\sym(\nabla u-P)$ {elastic strain}}}
\put(-20,183){\footnotesize{$\varepsilon_p=\sym P$  { micro-strain}}}
\put(-20,179){\footnotesize{ $\alpha=-\Curl P$ { micro-dislocation}}}
\put(0,177.5){\vector(0,-1){18}}
\put(-20,177.5){\vector(-1,-2){34.25}}
\put(22,178){\vector(1,-2){34.5}}
\put(0,138){\oval(46,42)}
\put(-19,155){\footnotesize{\bf A $(\dev,\dev)$-more relaxed }}
\put(-19,151){\footnotesize{\bf micromorphic dislocation}}
\put(-19,147){\footnotesize{12 dof $(u,P)$, well-posed$^*$}}
\put(-19,143){\footnotesize{$\sigma$ symmetric,\quad $\sigma={\mathbb{C}}.\, \varepsilon_e$}}
\put(-19,139){\footnotesize{isotropic: 6+2 parameters}}
\put(-19,135){\footnotesize{no coupling: 3+2 parameters}}
\put(-19,131){\footnotesize{\it constitutive variables:}}
\put(-19,127){\footnotesize{$\varepsilon_e=\sym(\nabla u-P)$ }}
\put(-19,123){\footnotesize{$\dev\varepsilon_p=\dev\sym P$  }}
\put(-19,119){\footnotesize{$\dev\alpha=-\dev\Curl P$ }}
\thicklines
\put(-53,88){\oval(45,42)}
\put(-71,106){\footnotesize{\bf A $(\cdot,\sym)$-more relaxed }}
\put(-71,102){\footnotesize{\bf micromorphic dislocation}}
\put(-71,98){\footnotesize{12 dof $(u,P)$, well-posed?}}
\put(-71,94){\footnotesize{$\sigma$ symmetric,\quad $\sigma={\mathbb{C}}.\, \varepsilon_e$}}
\put(-71,90){\footnotesize{isotropic: 6+2 parameters}}
\put(-71,86){\footnotesize{no coupling: 4+2 parameters}}
\put(-71,82){\footnotesize{\it constitutive variables:}}
\put(-71,78){\footnotesize{$\varepsilon_e=\sym(\nabla u-P)$ }}
\put(-71,74){\footnotesize{$\varepsilon_p=\sym P$  }}
\put(-71,70) {\footnotesize{$\sym\alpha=-\sym\Curl P$ }}
\put(54,88){\oval(47,42)}
\put(35,105){\footnotesize{\bf A $(\cdot,\dev\sym)$-more relaxed }}
\put(35,101){\footnotesize{\bf micromorphic dislocation}}
\put(35,97){\footnotesize{12 dof $(u,P)$, well-posed?}}
\put(35,93){\footnotesize{$\sigma$ symmetric,\quad $\sigma={\mathbb{C}}.\, \varepsilon_e$}}
\put(35,89){\footnotesize{isotropic: 6+1 parameters}}
\put(35,85){\footnotesize{no coupling: 4+1 parameters}}
\put(35,81){\footnotesize{\it constitutive variables:}}
\put(35,77){\footnotesize{$\varepsilon_e=\sym(\nabla u-P)$ }}
\put(35,73){\footnotesize{$\varepsilon_p=\sym P$  }}
\put(35,69) {\footnotesize{$\dev\sym\alpha=-\dev\sym\Curl P$ }}
\put(-53,67){\vector(0,-1){17.5}}
\put(-53,28){\oval(47,42)}
\put(-72,45){\footnotesize{\bf A $(\dev,\sym)$-more relaxed }}
\put(-72,41){\footnotesize{\bf micromorphic dislocation}}
\put(-72,37){\footnotesize{12 dof $(u,P)$, not well-posed!}}
\put(-72,33){\footnotesize{$\sigma$ symmetric,\quad $\sigma={\mathbb{C}}.\, \varepsilon_e$}}
\put(-72,29){\footnotesize{isotropic: 6+2 parameters}}
\put(-72,25){\footnotesize{no coupling: 3+2 parameters}}
\put(-72,21){\footnotesize{\it constitutive variables:}}
\put(-72,17){\footnotesize{$\varepsilon_e=\sym(\nabla u-P)$ }}
\put(-72,13){\footnotesize{$\dev\varepsilon_p=\dev\sym P$  }}
\put(-72,9) {\footnotesize{$\sym\alpha=-\sym\Curl P$ }}
\put(56,67){\vector(0,-1){17.5}}
\put(56,28){\oval(50,42)}
\put(35,45){\footnotesize{\bf A $(\dev,\dev\sym)$-more relaxed }}
\put(35,41){\footnotesize{\bf micromorphic dislocation}}
\put(35,37){\footnotesize{12 dof $(u,P)$, not well-posed!}}
\put(35,33){\footnotesize{$\sigma$ symmetric,\quad $\sigma={\mathbb{C}}.\, \varepsilon_e$}}
\put(35,29){\footnotesize{isotropic: 6+1 parameters}}
\put(35,25){\footnotesize{no coupling: 3+1 parameters}}
\put(35,21){\footnotesize{\it constitutive variables:}}
\put(35,17){\footnotesize{$\varepsilon_e=\sym(\nabla u-P)$ }}
\put(35,13){\footnotesize{$\dev\varepsilon_p=\dev\sym P$  }}
\put(35,9) {\footnotesize{$\dev\sym\alpha=-\dev\sym\Curl P$ }}
      \put(-29,28){\vector(1,0){60}}
\end{picture}
\end{center}
\caption{Relation between possible relaxed micromorphic models. The non-well-posedness follows from the results in \cite{BNPS1,BNPS2,BNPS3}. }\label{morerelaxfig}
\end{figure}
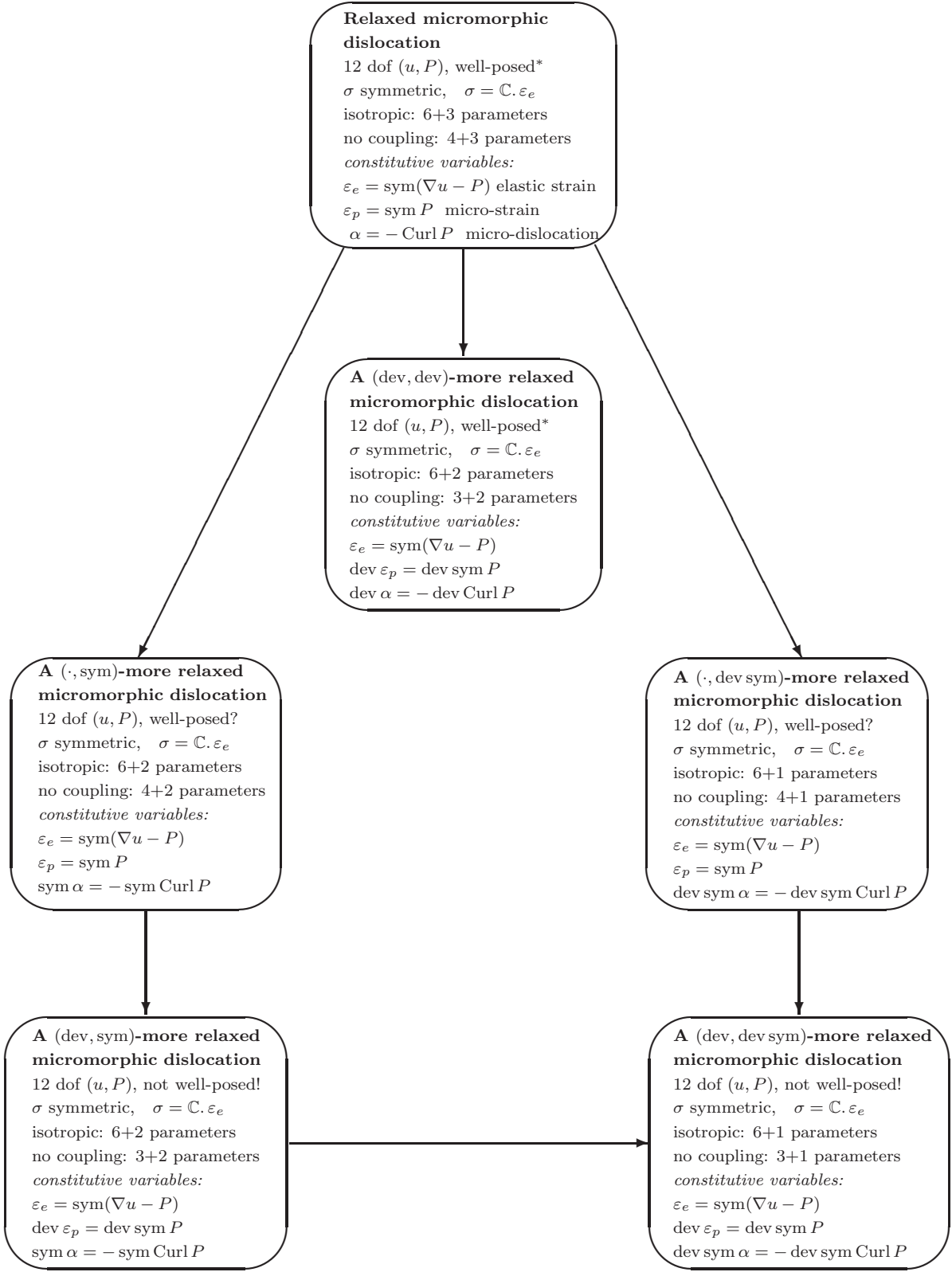

\section{Final remarks}
In the present paper we have mathematically studied a large class of evolution
equations which describe the behaviour of micromorphic
or generalized continua (see e.g. \cite{FdIPlacidiMadeo,PlacidiRosiMadeo,RosiMadeoGuyader}).
The mathematical existence, uniqueness and continuous dependence theorems
which we have obtained here are the logical basis of the studies which will
be developed in further investigations, where the manifold variety
of propagating mechanical waves which may exist in micromorphic continua
may unfold unexpected applications in the design of  particularly
tailored metamaterials \cite{Engheta,Zouhdi}, showing very useful and up-to-now unimagined
features. We  remark that the theorems obtained in the present paper
can also be used to give a better grounded basis to many results which
are already available in the literature (see e.g. \cite{FdIRosa,FdIRosa1}).

The diagram from Figure \ref{morerelaxfig} gives some new possible relaxed micromorphic models and, in view of   the status of the mathematical background, we indicate the well-posedness of the dynamic and static problem.

\bigskip

{\bf Acknowledgements.} I.D. Ghiba acknowledges support from the Romanian National Authority for Scientific Research (CNCS-UEFISCDI), Project No. PN-II-ID-PCE-2011-3-0521. I.D. Ghiba  would like to thank P. Neff  for his kind hospitality during his visit at the Faculty of Mathematics, Universit\"{a}t Duisburg-Essen, Campus Essen. P. Neff is grateful to F. dell'Isola for making his visit to CISTERNA di LATINA (M\&MoCS), in spring 2013, a wonderful scientific experience. A. Madeo thanks INSA-Lyon for the financial support assigned to the project BQR
2013-0054
``Mat\'{e}riaux M\'{e}so et Micro-H\'{e}t\'{e}rog\`{e}nes: Optimisation par Mod\`{e}les de Second Gradient et
Applications en Ing\'{e}nierie".

\bibliographystyle{plain} %plain
\addcontentsline{toc}{section}{References}

\begin{footnotesize}

\end{footnotesize}
\appendix
\section{Appendix}

\setcounter{equation}{0}
In this Appendix we outline some mathematical results used in this paper. We start with a generalization of Gronwall's Lemma:

\begin{lemma} {\rm [Brezis' Lemma (\cite[p. 47]{VrabieDiff}, Lemma 1.5.3)]}
 Let $x:[a,b]\rightarrow\mathbb{R}_+$ and $k:[a,b]\rightarrow\mathbb{R}_+$ be two continuous functions and let $m\geq 0$. If
 $
 \dd x^2(t)\leq m^2+2\int_0^tk(s)x(s) ds \ \ \text{for each}\ \  t\in[a,b],
 \ \ \text{then}$ \break $  \dd
 x(t)\leq m+\int_0^tk(s)ds  \ \ \text{for each}\ \  t\in[a,b]\, .
 $
 \end{lemma}

 In the following we present some important results regarding abstract Cauchy problems. Let us consider an Hilbert space $X$ and a linear operator $A$ with dense domain $\mathcal{D}(A)$ in $X$.

\begin{theorem} {\rm [Lumer--Phillips Corollary (\cite[p. 14]{Pazy}, Theorem 4.3)] }
 Let $A$ be a linear operator with dense domain $\mathcal{D}(A)$ in the Hilbert space $X$.
\begin{itemize}
\item[a)] If $A$ is dissipative and there is a $\lambda_0>0$ such that the range $R(\lambda_0I-A)$ of $\lambda_0I-A$ is $X$, then $A$ is the infinitesimal generator of a $C_0$ semigroup of contractions on $X$.
\item[b)] If $A$ is the infinitesimal generator of a $C_0$ semigroup of contractions on $X$, then  $R(\lambda I-A)=X$ for all $\lambda>0$ and $A$ is dissipative.
\end{itemize}
\end{theorem}
Let us consider a Hilbert space $X$ and the Cauchy Problem \cite{Pazy}
\begin{align}\label{PC}
\dd\frac{d}{dt}u=Au+f, \quad\quad\quad u(a)=\xi,
\tag{PC}
\end{align}
where $A:\mathcal{D}(A)\subseteq X\rightarrow X$ is the infinitesimal generator of a $C_0$ semigroup $\{S(t);t\geq 0\}$, $\xi\in X$ and $f\in L^1(a,b;X)$.

\begin{definition} {\rm (\cite[p. 183]{Vrabie})}
 The function $u:[a,b]\rightarrow X$ is called classical, or $C^1$-solution of the above problem \eqref{PC}, if $u$ is continuous on $[a,b]$, continuously differentiable on $(a,b]$, $u(t)\in \mathcal{D}(A)$ for each $t\in (a,b]$ and it satisfies $\dd\frac{d}{dt}u=Au+f$ for each $t\in [a,b]$ and $u(a)=\xi$.
\end{definition}

\begin{theorem}{\rm (\cite[p. 186]{Vrabie},  Corollary 8.1.2)}
 If $A:\mathcal{D}(A)\subseteq X\rightarrow X$ is the infinitesimal generator of a $C_0$ semigroup $\{S(t);t\geq 0\}$, and $f$ is of class $C^1$ on $[a,b]$, then, for each $\xi\in \mathcal{D}(A)$, the problem \eqref{PC} has a unique classical solution.
\end{theorem}

\begin{theorem}{\rm[Duhamel Principle (\cite[p. 184]{Vrabie}, Theorem 8.1.1)]}
 Each strong solution of \eqref{PC} is given by the variation of constants formula
\begin{equation}\label{vcf}
u(t)=S(t-a)\xi+\int_0^tS(t-s)f(s)ds.\tag{VCF}
\end{equation}
In particular, each classical solution of the problem \eqref{PC} is given by \eqref{vcf}.
\end{theorem}

\end{document}